\documentclass[11pt]{article}
\usepackage{amssymb}
\usepackage{amsmath}
\usepackage[all]{xy}
\usepackage{times}
\usepackage{graphicx}

\title{Analytic digraphs of uncountable Borel chromatic number\\ under injective definable homomorphism\indent}
\author{Dominique LECOMTE and Miroslav ZELEN\' Y$^1$}
\date{\today}

\def\ufootnote#1{\let\savedthfn\thefootnote\let\thefootnote\relax
\footnote{#1}\let\thefootnote\savedthfn\addtocounter{footnote}{-1}}

\newcommand{\Ana}{{\it\Sigma}^{1}_{1}}

\newcommand{\Ca}{{\it\Pi}^{1}_{1}}

\newcommand{\Borel}{{\it\Delta}^{1}_{1}}

\newcommand{\boraone}{{\bf\Sigma}^{0}_{1}}
\newcommand{\boratwo}{{\bf\Sigma}^{0}_{2}}

\newcommand{\borone}{{\bf\Delta}^{0}_{1}}

\newtheorem{thm} {Theorem} [section]
\newtheorem{defi} [thm] {Definition}
\newtheorem{cor} [thm] {Corollary}

\newtheorem{lem} [thm] {Lemma}

\begin{document}

\maketitle

\centerline{$\bullet$ Sorbonne Universit\'e, Institut de Math\'ematiques de Jussieu-Paris Rive Gauche,} 

\centerline{CNRS, Universit\'e Paris Diderot, Projet Analyse Fonctionnelle}

\centerline{Campus Pierre et Marie Curie, case 247, 4, place Jussieu, 75 252 Paris cedex 5, France}

\centerline{dominique.lecomte@upmc.fr}\bigskip

\centerline{$\bullet$ Universit\'e de Picardie, I.U.T. de l'Oise, site de Creil,}

\centerline{13, all\'ee de la fa\"\i encerie, 60 107 Creil, France}\bigskip

\centerline{$\bullet^1$ Charles University, Faculty of Mathematics and Physics, Department of Mathematical Analysis}

\centerline{Sokolovsk\'a 83, 186 75 Prague, Czech Republic}

\centerline{zeleny@karlin.mff.cuni.cz}\bigskip\bigskip\bigskip\bigskip\bigskip\bigskip

\ufootnote{{\it 2010 Mathematics Subject Classification.}~03E15, 54H05}

\ufootnote{{\it Keywords and phrases.}~analytic, antichain, Borel, chromatic number, digraph, homomorphism, injective, minimal}

\ufootnote{${}^1$ The second author was supported by the grant GA\v{C}R 15-08218S.}

\noindent {\bf Abstract.} We study the analytic digraphs of uncountable Borel chromatic number on Polish spaces, and compare them with the notion of injective Borel homomorphism. We provide some minimal digraphs incomparable with $\mathbb{G}_0$. We also prove the existence of antichains of size continuum, and that there is no finite basis.

\vfill\eject

\section{$\!\!\!\!\!\!$ Introduction}\indent

 In this paper, $A,B$ will be (binary) relations on some sets $X,Y$ respectively. The {\bf diagonal} of $X$ is $\Delta (X)\! :=\!\{ (x,x)\mid x\!\in\! X\}$. We will say that $A$ is {\bf irreflexive}, or a {\bf digraph}, if $A$ does not meet the diagonal (some authors call these relations {\bf simple digraphs}). We set 
${A^{-1}\! :=\!\{ (x,y)\!\in\! X^2\mid (y,x)\!\in\! A\}}$. We say that $A$ is {\bf symmetric} if $A\! =\! A^{-1}$, and $A$ is {\bf antisymmetric} if 
$A\cap A^{-1}\!\subseteq\!\Delta (X)$. The set $s(A)\! :=\! A\cup A^{-1}$ is the {\bf symmetrization} of $A$. We say that $A$ is a {\bf graph} if $A$ is irreflexive and symmetric, and $A$ is an {\bf oriented graph} if $A$ is irreflexive and antisymmetric. Recall that $A$ is a {\bf subgraph} of the digraph $B$ if 
$X\!\subseteq\! Y$ and $A\!\subseteq\! B$. An $A$-{\bf path} is a finite sequence $(x_i)_{i\leq n}$ of points of $X$ such that $(x_i,x_{i+1})\!\in\! A$ if $i\! <\! n$. We say that $A$ is {\bf connected} if for any $x,y\!\in\! X$ there is an $A$-path $(x_i)_{i\leq n}$ with $x_0\! =\! x$ and $x_n\! =\! y$. If $A$ is a graph, then the 
{\bf connected component} of $x\!\in\! X$ is the set $C(x)$ of points $y\!\in\! X$ for which there is an $A$-path $(x_i)_{i\leq n}$ with $x_0\! =\! x$ and 
$x_n\! =\! y$. A graph $A$ is {\bf acyclic} if there is no injective $A$-path $(x_i)_{i\leq n}$ with $n\!\geq\! 2$ and $(x_n,x_0)\!\in\! A$. If $A$ is an acyclic graph and $y\!\in\! C(x)$, then we denote by $p^A_{x,y}$ the unique injective $A$-path $(x_i)_{i\leq n}$ with $x_0\! =\! x$ and $x_n\! =\! y$. We say that $A$ is {\bf locally countable} if its horizontal and vertical sections are countable.\bigskip

 We write $(X,A)\preceq (Y,B)$ when there is $h\!:\! X\!\rightarrow\! Y$ such that 
$A\!\subseteq\! (h\!\times\! h)^{-1}(B)$. If this holds, then we say that $h$ is a 
{\bf homomorphism} from $(X,A)$ into $(Y,B)$. When $h$ can be injective, we write $(X,A)\preceq^{\mbox{inj}}(Y,B)$. The notion of injective homomorphism is very natural since it corresponds to the basic notion of subgraph. Indeed, if $h$ is an injective homomorphism from $(X,A)$ into $(Y,B)$, then 
$(h[X],(h\!\times\! h)[A])$ is a subgraph of $(Y,B)$. Conversely, if $(X,A)$ is a subgraph of $(Y,B)$, then the canonical injection is an injective homomorphism from $(X,A)$ into $(Y,B)$. A {\bf coloring} from $(X,A)$ into some set $Y$ is  a map $c\! :\! X\!\rightarrow\! Y$ such that $c(x)\!\not=\! c(x')$ if $(x,x')\!\in\! A$, i.e., a homomorphism from $(X,A)$ into $(Y,\not= )$.\bigskip

 The reader should see [K] for the standard descriptive set theoretic notions and notation. Let $\mathfrak C$ be a class of functions between Polish spaces, e.g., continuous, Borel (denoted $c,B$ respectively). If $X,Y$ are Polish spaces and $h$ can be in $\mathfrak C$, then we will use the notation $\preceq^{}_{\mathfrak C}$, 
$\preceq^{\mbox{inj}}_{\mathfrak C}$ respectively. The study of definable colorings of analytic graphs was initiated in [K-S-T]. The $\mathfrak C$-{\bf chromatic number} of a digraph $A$ on a Polish space $X$ is the smallest cardinality of a Polish space $Y$ for which there is a $\mathfrak C$-coloring from $(X,A)$ into $Y$.\bigskip
 
\noindent\bf Example.\rm\ Let $\psi\! :\!\omega\!\rightarrow\! 2^{<\omega}$ be a natural bijection ($\psi (0)\! =\!\emptyset$, $\psi (1)\! =\! 0$, $\psi (2)\! =\! 1$, 
$\psi (3)\! =\! 0^2$, $\psi (4)\! =\! 01$, $\psi (5)\! =\! 10$, $\psi (6)\! =\! 1^2$, $\ldots$). A crucial property of $\psi$ is that $\psi^{-1}(s)\! <\!\psi^{-1}(s\varepsilon )$ if $s\!\in\! 2^{<\omega}$ and $\varepsilon\!\in\! 2$.\bigskip

 Note that $|\psi (n)|\!\leq\! n$, so that we can define $s_n\! :=\!\psi (n)0^{n-|\psi (n)|}$. Some crucial properties of 
$(s_n)$ are that it is {\bf dense} (for each $s\!\in\! 2^{<\omega}$, there is $n$ such that $s\!\subseteq\! s_n$), and that $|s_{n}|\! =\! n$. We set 
$\mathbb{G}_0\!:=\!\{ (s_{n}0\gamma ,s_{n}1\gamma )\mid n\!\in\!\omega\wedge\gamma\!\in\! 2^\omega\}$. The set $s(\mathbb{G}_0)$ is considered in [K-S-T], where the following is essentially proved.

\begin{thm} (Kechris, Solecki, Todor\v cevi\' c) \label{G0} Let $X$ be a Polish space and $A$ be an analytic digraph on $X$. Then exactly one of the following holds:\smallskip

(a) $(X,A)\preceq_{B}(\omega ,\not=)$ (i.e., $(X,A)$ has countable Borel chromatic number),\smallskip

(b) $(2^\omega ,\mathbb{G}_0)\preceq_{c}(X,A)$.\end{thm}

 Actually, the original statement in [K-S-T] is when $A$ is a graph, with $s(\mathbb{G}_0)$ instead of 
$\mathbb{G}_0$. But we can get Theorem \ref{G0} without any change in the proof in [K-S-T].\bigskip

  This result had a lot of developments since. For instance, Miller developed some techniques to recover many dichotomy results of descriptive set theory, without using effective descriptive set theory (see [Mi]). He replaces it with some versions of Theorem \ref{G0}. In [K-S-T], it is conjectured that we can replace $\preceq_c$ with $\preceq^{\mbox{inj}}_c$ in the version of Theorem \ref{G0}.(b) for graphs (the authors show in Theorem 6.6  that this is the case if $A$ is an acyclic graph or a locally countable graph, just like $s(\mathbb{G}_0)$; their proof also works for digraphs with acyclic symmetrization or locally countable, with $\mathbb{G}_0$ instead of $s(\mathbb{G}_0)$). It is proved in [L4] that this is not the case.

\begin{thm} \label{nikst} (Lecomte) There is no pair $(\mathbb{X}_0,\mathbb{A}_0)$, where 
$\mathbb{X}_0$ is Polish and $\mathbb{A}_0$ is an analytic graph on $\mathbb{X}_0$, such that for every pair $(X,A)$ of the same type, exactly one of the following holds:\smallskip

(a)~$(X,A)\preceq_{B}(\omega ,\not= )$,\smallskip

(b)~$(\mathbb{X}_0,\mathbb{A}_0)\preceq^{\mbox{inj}}_c(X,A)$.\end{thm}

 In other words, there is no one-element basis for $\preceq^{\mbox{inj}}_c$ among analytic graphs of uncountable Borel chromatic number (recall that if 
$(Q,\leq )$ is a quasi-ordered space, then a {\bf basis} is a subfamily $\cal F$ of $Q$ such that any element of $Q$ is $\leq$-above an element of $\cal F$). This led Kechris and Marks to ask the following in [K-Ma] (see Problem 3.39).\bigskip

\noindent\bf Questions\rm\ (1) Is there a basis of cardinality $<\! 2^{\aleph_0}$ for $\preceq^{\mbox{inj}}_B$ among analytic graphs of uncountable Borel chromatic number?\smallskip

\noindent (2) If not, is there such a basis consisting of a continuum size family of Òreasonably simpleÓ graphs?\bigskip

 Of course, we can ask the same questions with $\preceq^{\mbox{inj}}_c$ instead of 
$\preceq^{\mbox{inj}}_B$, and for digraphs instead of graphs. We are interested in basis as small as possible with respect to the inclusion. In other words, we want our basis to be antichains for the quasi-order we consider (recall that a subfamily $\cal F$ of $Q$ is an {\bf antichain} if the elements of $\cal F$ are pairwise $\leq$-incomparable). This leads to the following.\bigskip

\noindent\bf Question\rm\ (3) Is there a $\preceq^{\mbox{inj}}_{\mathfrak C}$-antichain basis for the class of analytic digraphs of uncountable Borel chromatic number?\bigskip

 In [L-Mi], it is proved that there is neither $\leq_{\mathfrak C}$-antichain basis, nor $\sqsubseteq_{\mathfrak C}$-antichain basis, for the class of analytic graphs of uncountable Borel chromatic number (these two quasi-orders are defined like $\preceq_{\mathfrak C}$ and $\preceq^{\mbox{inj}}_{\mathfrak C}$ respectively, except that 
``$A\!\subseteq\! (h\!\times\! h)^{-1}(B)$" is replaced with ``$A\! =\! (h\!\times\! h)^{-1}(B)$"). In particular, (3) may have a negative answer. A first approach for the result in [L-Mi] was the existence of an $\leq_{\mathfrak C}$-antichain of size $2^{\aleph_0}$ made of graphs $\leq_{\mathfrak C}$-minimal among analytic graphs of uncountable Borel chromatic number, essentially proved in [L4]. This leads to the following.\bigskip

\noindent\bf Question\rm\ (4) Is there a $\preceq^{\mbox{inj}}_{\mathfrak C}$-antichain of size $2^{\aleph_0}$ made of digraphs $\preceq^{\mbox{inj}}_{\mathfrak C}$-minimal among analytic digraphs of uncountable Borel chromatic number?\bigskip

 The minimal elements of $(Q,\leq )$ are of particular importance since they have to be part of any basis, up to equivalence. The discussion after Theorem \ref{G0} shows that $\mathbb{G}_0$ is $\preceq_c^{\mbox{inj}}$ and 
$\preceq_B^{\mbox{inj}}$-minimal among analytic digraphs of uncountable Borel chromatic number. The main results in this paper are steps towards a positive answer to Question (4), and our paper is mentioned in the last version of [K-Ma]. 
 
\begin{thm} \label{G1} Let ${\mathfrak C}\!\in\!\{ c,B\}$. There is a $\preceq^{\mbox{inj}}_{\mathfrak C}$-antichain 
$\{\mathbb{G}_0,\mathbb{G}_1,\mathbb{G}_1^{-1}\}$ made of digraphs $\preceq^{\mbox{inj}}_{\mathfrak C}$-minimal among analytic digraphs of uncountable Borel chromatic number.\end{thm}

 Thus Theorem \ref{nikst} holds for digraphs. Note that some antichains were already present in [L4]. The main result in the present paper is the minimality of $\mathbb{G}_1$, and thus a dichotomy result. We now provide a construction of the digraph $\mathbb{G}_1$. Note that $\mathbb{G}_0$ can be viewed as the union of the graphs of the partial homeomorphisms $h_n\! :\! N_{s_n0}\!\rightarrow\! N_{s_n1}$ given by 
$h_n(s_n0\gamma )\! :=\! s_n1\gamma$. In particular, the following picture holds if $m\! <\! n$ and 
$h_m\big( h_n(\alpha )\big)$, $h_m(\alpha )$ are defined:
$$\xymatrix{ 
\alpha\ar[d]_{h_n}\ar[r]^{h_m~~~~} & h_m(\alpha )\ar@{}[d]_{\not=}\\ 
h_n(\alpha )\ar[r]^{h_m~~~~~~}  & h_m\big( h_n(\alpha )\big)
}$$
We define $t_n\!\in\! 2^{<\omega}$ and maps $g_n\! :\! N_{t_n0}\!\rightarrow\! N_{t_n1}$, and $\mathbb{G}_1$ will be the union of the graphs of the $g_n$'s. One of the crucial properties of the $g_n$'s is that 
$g_m\big( g_n(\alpha )\big)\! =\! g_m(\alpha )$ if $m\! <\! n$ and $g_m\big( g_n(\alpha )\big)$, $g_m(\alpha )$ are defined. In particular, the following picture holds and violates the previous one:
$$\xymatrix{ 
\alpha\ar[d]_{g_n}\ar[dr]^{g_m~~~~}\\ 
g_n(\alpha )\ar[r]^{g_m~~~~~~~~~~~~~}  & g_m\big( g_n(\alpha )\big)\! =\! g_m(\alpha )
}$$
This provides some cycles, which have to exist in examples orthogonal to $\mathbb{G}_0$, by the discussion after Theorem \ref{G0}. The $\leq_{\mathfrak C}$-antichain mentioned before Question (4) was constructed with different configurations of cycles. We believe that some other algebraic conditions of this type could lead to a positive answer to Question (4). Our second main result is a weak version of this.

\begin{thm} \label{biganti} Let ${\mathfrak C}\!\in\!\{ c,B\}$.\smallskip

(a) There is a $\preceq^{\mbox{inj}}_{\mathfrak{C}}$-antichain of size $2^{\aleph_0}$ made of 
$\boratwo$ digraphs of uncountable Borel chromatic number.\smallskip

(b) There is a $\preceq^{\mbox{inj}}_{\mathfrak{C}}$-strictly increasing chain of size $\aleph_0$ made of 
$\boratwo$ digraphs of uncountable Borel chromatic number.\end{thm}

 These digraphs, in fact differences of two closed sets, are not all minimal. In fact, we first construct a 
$\preceq^{\mbox{inj}}_B$-antichain of size $\aleph_0$ made of digraphs in the style of $\mathbb{G}_1$, that could be minimal as in Theorem \ref{G1}. We then consider suitable direct sums of these digraphs. This antichain of size $\aleph_0$ is in fact made of pairwise incompatible digraphs, which gives our third main result (recall that 
$p,q\!\in\! Q$ are {\bf incompatible} if there is no $r\!\in\! Q$ with $r\!\leq\! p,q$).

\begin{thm} \label{infi} Let ${\mathfrak C}\!\in\!\{ c,B\}$. Any $\preceq^{\mbox{inj}}_{\mathfrak{C}}$-basis for the class of analytic digraphs of uncountable Borel chromatic number on Polish spaces is infinite.\end{thm}

 We now provide a concrete description of $\mathbb{G}_1$. We first define a sequence $(q_n)_{n\in\omega}$ of natural numbers by setting $q_0\! :=\! 0$ and $q_{n+1}\! :=\! 32^{q_n}$. Note that $(q_n)_{n\in\omega}$ is strictly increasing. In particular, $\vert\psi (n)\vert\!\leq\! n\!\leq\! q_n\! <\! 2^{q_n}$, so that 
$t_n\! :=\!\psi (n)0^{2^{q_n}-\vert\psi (n)\vert}$ is well-defined and has length $2^{q_n}$. We then set 
$S_n\! :=\!\{ 2^{q_n}\!\cdot\! j\mid j\!\geq\! 1\}$, define 
$\theta_n\! :\!\omega\!\rightarrow\!\omega$ and $g_n$ by 
$$\theta_n(k)\! :=\!\left\{\!\!\!\!\!\!\!\!
\begin{array}{ll} 
& k\mbox{ if }k\!\notin\! S_n\mbox{,}\cr
& 2^{q_n}\!\cdot\! (2j\! +\! 1)\mbox{ if }k\!\ =\! 2^{q_n}\!\cdot\! j\mbox{,}
\end{array}
\right.
~~~~~~~~~~
g_n(\alpha )(k)\! :=\!\left\{\!\!\!\!\!\!\!\!
\begin{array}{ll} 
& 1\mbox{ if }k\! =\! 2^{q_n}\mbox{,}\cr 
& \alpha\big(\theta_n(k)\big)\mbox{ if }k\!\not=\! 2^{q_n}.
\end{array}
\right.$$
The next definition catches some of the crucial properties of the sequence $(g_n)$ defining $\mathbb{G}_1$. As noted in [L4], a great variety of very different non-potentially closed relations appear at the level of differences of two closed sets. For the kind of examples we will consider, being non-potentially closed and having uncountable Borel chromatic number are equivalent properties. Let us make this more precise.\bigskip
 
\noindent\bf Notation.\rm\ If $(f_n)$ is a sequence of functions, then we set 
$A^f\! :=\!\bigcup_{n\in\omega}~\mbox{Graph}(f_{n})$.

\begin{defi}  We say that $\big( X,(f_n)\big)$ is a {\bf complex situation} if\smallskip

(a) $X$ is a nonempty Polish space,\smallskip
  
(b) the $f_n$'s are partial continuous and open maps whose domain and range are open in $X$,\smallskip

(c) $\Delta (X)\!\subseteq\!\overline{A^f}\!\setminus\! A^f$.\end{defi}

 This kind of situations play an important role in the theory of potential complexity (see, for example, Definition 2.2 in [L3], and also Definitions 13, 26 and 31 in [L4]). These properties are sufficient to ensure that $A^f$ is a $\boratwo$ digraph of uncountable Borel chromatic number (see Corollary \ref{D2+}). We prove a  result giving some additional motivation for introducing this notion, which is in fact very general.
 
\begin{thm} \label{fav} Let $Y$ be a Polish space and $B$ be an analytic digraph on $Y$. Then exactly one of the following holds:\smallskip

(a) $(Y,B)$ has countable Borel chromatic number,\smallskip

(b) there is a complex situation $\big( X,(f_n)\big)$ such that the inequality $(X,A^f)\preceq_c^{\mbox{inj}}(Y,B)$ or 
$(X,A^f)\preceq_c^{\mbox{inj}}(Y,B^{-1})$ holds.\end{thm}

 We prove more than the minimality of $\mathbb{G}_1$.

\begin{thm} \label{cormain} Let $\big( X,(f_n)\big)$ be a complex situation satisfying the following additional property:\smallskip

(d) $f_m\big( f_n(x)\big)\! =\! f_m(x)$ if $m\! <\! n$ and $f_m\big( f_n(x)\big) ,f_m(x)$ are defined.\smallskip

\noindent Then $(2^\omega ,\mathbb{G}_1)\preceq^{\mbox{inj}}_c(X,A^f)$.\end{thm}

 The organization of this paper is as follows. In Section 2, we provide some basic properties of complex situations. In Section 3, we characterize when the digraph associated with a complex situation is minimal among analytic digraphs of uncountable Borel chromatic number, and prove Theorem \ref{fav}. In Section 4, we prove a relatively general lemma ensuring the injectivity of the homomorphism $h$ implicitly mentioned in the statement of Theorem \ref{cormain}.
 
\vfill\eject
 
 In Section 5, we show that $\mathbb{G}_1$ comes from a complex situation, introduce some finitary objects used in the construction of $h$, and prove their important properties. In particular, in many Cantor-like constructions of homomorphisms or reductions, we construct approximations, indexed by finite binary sequences, of the desired infinitary objects. The construction is usually made by induction on the length of the finite binary sequences. So we consider the partitions into basic clopen sets $(N_x)_{x\in 2^l}$ of $2^\omega$, for each $l\!\in\!\omega$. Here it will be more convenient to replace $2^l$ with some subset $X_l$ of $2^{<\omega}$ containing sequences of different lengths since the $g_n$'s ``forget" some coordinates. In Section 6, we construct our homomorphism and prove Theorems \ref{G1} and \ref{cormain}. In Section 7, we prove Theorems \ref{biganti} and \ref{infi}, and show that our three main results also hold for graphs.

\section{$\!\!\!\!\!\!$ Some basic properties of complex situations}\indent

 The next lemma is essentially Lemma 3.5 in [L1], and the crucial point of its proof.

\begin{lem} \label{suffnpc} Let $X$ be a nonempty Polish space and, for $n\!\in\!\omega$, $D_n,R_n$ be dense $G_{\delta}$ subsets of some open subsets of $X$, and $f_n\! :\! D_n\!\rightarrow\! R_n$ be a continuous, open and onto map.\smallskip

(a) Let $G$ be a dense $G_{\delta}$ subset of $X$. Then $\mbox{Graph}(f_{n})\!\subseteq\!\overline{\mbox{Graph}(f_{n})\cap G^{2}}$.\smallskip

(b) We assume that $\Delta (X)\!\subseteq\!\overline{A^f}\!\setminus\! A^f$. Then $(X,A^f)$ has uncountable Borel chromatic number.\end{lem}

\noindent\bf Proof.\rm\ (a) See Lemma 3.5 in [L1].\bigskip

\noindent (b) As $\Delta (X)\!\subseteq\!\neg A^f$, $(X,A^f)$ is a digraph and thus its Borel chromatic number is defined. We argue by contradiction to see that it is uncountable, which gives a countable partition $(B_n)$ of $X$ into Borel sets with $A^f\cap B_n^2\! =\!\emptyset$. By 13.5 in [K], there is a finer Polish topology $\tau$ on $X$ such that the $B_n$'s are clopen in $(X,\tau )$. By 15.2 in [K], the identity of $X$ equipped with its initial topology into $(X,\tau )$ is Borel. By 11.5 in [K] it is Baire measurable. By 8.38, in [K], there is a dense $G_\delta$ subset $G$ of $X$ on which the $B_n$'s are clopen. We pick $x\!\in\! G$, which exists since $X$ is Polish and nonempty. We choose $n$ with $x\!\in\! B_n$. By (a), 
$\overline{A^f}\cap G^2\! =\!\overline{A^f\cap G^2}\cap G^2$, so that 
$(x,x)\!\in\!\overline{A^f\cap G^2}\cap (B_n\cap G)^2$. As $B_n\cap G$ is clopen in $G$, 
$A^f\cap (B_n\cap G)^2$ is nonempty, which is absurd.\hfill{$\square$}

\begin{cor} \label{D2+} Let $\big(X,(f_{n})\big)$ be a complex situation. Then $A^f$ is a $\boratwo$ digraph of uncountable Borel chromatic number.\end{cor}

\noindent {\bf Proof.}~Note that $A^f\! =\!\bigcup_{n\in\omega}\mbox{ Graph}(f_n)$ is $\boratwo$ since the $f_n$'s are continuous with open domain. Moreover, Lemma \ref{suffnpc} ensures that $(X,A^f)$ has uncountable Borel chromatic number.\hfill{$\square$}\bigskip

 Condition (d) in Theorem \ref{cormain} is motivated by the following result, which is essentially Claim 1 in the proof of Theorem 10 in [L4]. This condition ensures that $A^f$ is $\preceq^{\mbox{inj}}_B$-incomparable with 
$\mathbb{G}_0$.\bigskip
 
\noindent\bf Notation.\rm\ In the sequel, $D_n$ will be the domain of the function $f_n$.

\begin{lem} \label{count} Let $X$ be a Polish space, and $g_0,g_1,\ldots\! :\! X\!\rightarrow\! X$ be fixed point free Borel partial functions such that $g_m\big( g_n(x)\big)\! =\! g_m(x)$ if $m\! <\! n$ and $g_m\big( g_n(x)\big) ,g_m(x)$ are defined. Then every locally countable analytic subset of $A^g$ has countable Borel chromatic number.\end{lem} 

\noindent {\bf Proof.} Suppose that $H$ is a locally countable analytic subset of $A^g$. By 35.13 in [K] and Lemma 2.4.(a) in [L2], there are Borel partial injections $f_n$ on $X$ such that 
$H\!\subseteq\! A^f\!\subseteq\! A^g$. By replacing each $f_n$ with its restrictions to the sets 
$\{ x\in D_n\mid f_n(x)\! =\! g_m(x)\}$, for $m\in\omega$, we can assume that, for all $n\in\omega$, there is $k_n\!\in\!\omega$ such that $f_n\! =\! {g_{k_n}}_{\vert_{D_n}}$. It is easily seen that the graph of a fixed point free Borel function has countable Borel chromatic number (see Proposition 4.5 of [K-S-T]). So, by replacing $f_n$ with its restriction to countably many Borel sets, we can also assume that 
$D_n^2\cap\bigcup_{k\leq k_n}\mbox{ Graph}(g_k)\! =\! \emptyset$ for all $n\in\omega$. It only remains to note that $D_n^2\cap\bigcup_{k>k_n}\mbox{ Graph}(g_k)=\emptyset$. In order to see this, simply observe that if 
$k\! >\! k_n$ and $x,g_k(x)\in D_n$, then $f_n(x)\! =\! g_{k_n}(x)\! =\! g_{k_n}\circ g_k(x)\! =\! f_n\circ g_k(x)$, which contradicts the fact that $f_n$ is a partial injection.\hfill{$\square$}

\begin{cor} \label{th12} Let $\big( X,(f_n)\big)$ be a complex situation satisfying Condition (d) in Theorem \ref{cormain}. Then $\{\mathbb{G}_0,A^f,(A^f)^{-1}\}$ is a $\preceq^{\mbox{inj}}_B$-antichain.\end{cor}

\noindent {\bf Proof.} Assume that $\mathbb{G}_0$ is $\preceq^{\mbox{inj}}_B$-below $A^f$, with witness $\pi$. Then $(\pi\times\pi )[\mathbb{G}_0]$ is a locally countable Borel subset of $A^f$ with uncountable Borel chromatic number, which contradicts Lemma \ref{count}. $A^f$ is not $\preceq^{\mbox{inj}}_B$-below $\mathbb{G}_0$ since $\mathbb{G}_0$ is locally countable and $A^f$ is not, by Corollary \ref{D2+} and Lemma \ref{count}. The discussion after Theorem \ref{G0} shows that 
$(2^\omega ,\mathbb{G}_0)\preceq^{\mbox{inj}}_c(2^\omega ,\mathbb{G}_0^{-1})$. Thus 
$\mathbb{G}_0$ and $(A^f)^{-1}$ are $\preceq^{\mbox{inj}}_B$-incomparable. If 
$(2^\omega ,(A^f)^{-1})\preceq^{\mbox{inj}}_B(2^\omega ,A^f)$ with witness $h$, then 
$(h\!\times\! h)[(A^f)^{-1}]$ is a locally countable subset of $A^f$ since $A^f$ has countable vertical sections, which contradicts Lemma \ref{count}.\hfill{$\square$}\bigskip

 The next two lemmas will be used in the proof of Theorem \ref{cormain}.
 
\begin{lem} \label{trans} Let $\big( X,(f_n)\big)$ be a complex situation satisfying Condition (d) in Theorem \ref{cormain}, $V_0,V_1$ be subsets of $X$, and $m\! <\! n$ be natural numbers such that 
$V_0\!\subseteq\! D_m\cap D_n$ and $V_1\!\subseteq\! f_n[V_0]\cap D_m$. Then 
$f_m[V_1]\!\subseteq\! f_m[V_0]$.\end{lem}

\noindent\bf Proof.\rm\ Pick $y\!\in\! V_1$, and $x\!\in\! V_0$ with $y\! =\! f_n(x)$. Note that $f_m\big( f_n(x)\big)$ is defined, as well as $f_m(x)$. By Condition (d), $f_m\big( f_n(x)\big)\! =\! f_m(x)$. This implies that 
$f_m(y)\! =\! f_m(x)\!\in\! f_m[V_0]$.\hfill{$\square$}

\begin{lem} \label{creation+} Let $\big( X,(f_n)\big)$ be a complex situation such that $X$ is zero-dimensional and the $D_n$'s are clopen, $V$ be a nonempty open subset of $X$, and $m$ be a natural number. Then we can find $n\! >\! m$ and nonempty clopen subsets $V_0,V_1$ of $X$ such that 
$V_0\!\subseteq\! V\cap D_n$ and $V_1\!\subseteq\! V\cap f_n[V_0]$.\end{lem}

\noindent\bf Proof.\rm\ The assumption on $\big( X,(f_n)\big)$ implies that $\Delta (X)\!\subseteq\!\overline{\bigcup_{n>m}~\mbox{Graph}(f_n)}$ since the 
$\mbox{Graph}(f_n)$'s are closed. This gives $n\! >\! m$ such that $V^2\cap\mbox{Graph}(f_n)\!\not=\!\emptyset$, and $(x,y)$ in this intersection. In particular, $x\!\in\! D_n$ and $y\! =\! f_n(x)$. We choose a clopen subset $V_0$ of $X$ with $x\!\in\! V_0\!\subseteq\! V\cap D_n$, and a clopen subset $V_1$ of $X$ with 
$y\!\in\! V_1\!\subseteq\! V\cap f_n[V_0]$.\hfill{$\square$}

\section{$\!\!\!\!\!\!$ The characterization of the minimality}\indent

 We will characterize when the set $A^f$ associated with a complex situation $\big( X,(f_n)\big)$ is 
$\preceq^{\mbox{inj}}_c$-minimal among analytic digraphs of uncountable Borel chromatic number. We will need a strengthening of the notion of a complex situation.

\begin{defi}  We say that a complex situation $\big( X,(f_n)\big)$ is a {\bf strongly complex situation} if\smallskip

(a) $X$ is a nonempty zero-dimensional perfect Polish space,\smallskip
  
(b) the $f_n$'s are partial continuous and open maps whose domain and range are clopen in $X$,\smallskip
  
(c) the restriction of any $f_n$ to any nonempty open subset of its domain is not countable-to-one.\end{defi}

 The reader should see [M] for the basic notions of effective descriptive set theory. Let $X$ be a recursively presented Polish space. The topology 
${\it\Delta}_X$ on $X$  is generated by $\Borel (X)$. This topology is Polish (see the proof of Theorem 3.4 in [Lo]). The Gandy-Harrington topology 
${\it\Sigma}_X$ on $X$ is generated by $\Ana (X)$. Recall that $\Omega_X\! :=\!\{ x\!\in\! X\mid\omega_1^x\! =\!\omega_1^{\mbox{CK}}\}$ is Borel and $\Ana$, and $(\Omega_X,{\it\Sigma}_X)$ is a zero-dimensional Polish space (in fact, the intersection of $\Omega_X$ with any nonempty $\Ana$ set is a nonempty clopen subset of $(\Omega_X,{\it\Sigma}_X)$-see [L1]).

\begin{lem} \label{Bas0+} Let $\big(X,(f_n)\big)$ be a complex situation such that 
$(2^\omega ,\mathbb{G}_0)\not\preceq_c^{\mbox{inj}}(X,A^f)$, $P$ be a Borel subset of $X$, and $(S_n)$ be a sequence of analytic subsets of $X$ such that 
$$A\! :=\!\bigcup_{n\in\omega}~\mbox{Graph}({f_n}_{\vert S_n})\!\subseteq\! P^2$$ 
has uncountable Borel chromatic number. Then we can find a Borel subset $S$ of $P$, a finer topology $\tau$ on 
$S$, and a sequence $(C_n)$ of clopen subsets of $Y\! :=\! (S,\tau )$ such that $\big( Y,({f_n}_{\vert C_n})\big)$ is a strongly complex situation and $C_n\!\subseteq\! S\cap f_n^{-1}(S)\cap S_n$ for each $n\!\in\!\omega$.\end{lem}

\noindent {\bf Proof.} In order to simplify the notation, we will assume as usual that $X$ is recursively presented, 
$(f_n)$ is $\Borel$ (so that $A^f$ is $\Borel$ too), $P$ is $\Borel$, and $(S_n)$ is $\Ana$ (so that $A$ is $\Ana$). We set, for each $n\!\in\!\omega$ and each $W\!\subseteq\! X$, 
$$\Phi_n(W)\Leftrightarrow {f_n}_{\vert D_n\cap W}\mbox{ is countable-to-one.}$$ 
Note that $\Phi_n$ is $\Ca$ on $\Ana$. Indeed, let $Z$ be a recursively presented Polish space, and 
$\mathfrak S$ be in $\Ana (Z\!\times\! X)$. Then $\Phi_n({\mathfrak S}_z)\Leftrightarrow
\forall x\!\in\! X~~f_n^{-1}(\{ x\} )\cap {\mathfrak S}_z\mbox{ is countable}$. Note that 
$f_n^{-1}(\{ x\} )\cap {\mathfrak S}_z$ is $\Ana (z,x)$. By 4F.1 in [M], $f_n^{-1}(\{ x\} )\cap {\mathfrak S}_z$ is countable if and only if it is contained in $\Borel (z,x)\cap X$, which is a $\Ca$ condition (in $(z,x)$).\bigskip

 This argument shows that if $X_n\! :=\!\bigcup\{\Delta\!\in\!\Borel (X)\mid\Phi_n(\Delta)\}$, then $(X_n)$ is 
$\Ca$, and also that we can apply the effective version of the first reflection theorem (see 35.10 in [K]). Let us prove that if $C\!\in\!\Ana (X)$ and $C\!\setminus\! X_n$ is not empty, then ${f_n}_{\vert D_n\cap C}$ is not countable-to-one. We argue by contradiction. As $C\!\setminus\! X_n\!\in\!\Ana$ and 
$\Phi_n(C\!\setminus\! X_n)$ holds, the effective version of the first reflection theorem gives $\Delta\!\in\!\Borel$ such that $C\!\setminus\! X_n\!\subseteq\!\Delta$ and $\Phi_n(\Delta )$ holds. Thus $\Delta\!\subseteq\! X_n$ and $C\!\setminus\! X_n\!\subseteq\! X_n\!\setminus\! X_n$ is empty, which is absurd.\bigskip

 Note that $A\! =\!\bigcup_{n\in\omega}~\mbox{Graph}({f_n}_{\vert S_n\cap X_n})\cup
\bigcup_{n\in\omega}~\mbox{Graph}({f_n}_{\vert S_n\setminus X_n})$. As 
$\bigcup_{n\in\omega}~\mbox{Graph}({f_n}_{\vert S_n\cap X_n})$ is locally countable analytic, it has countable Borel chromatic number since $(2^\omega ,\mathbb{G}_0)\not\preceq_c^{\mbox{inj}}(X,A^f)$, by the discussion after Theorem \ref{G0}. Thus $A'\! :=\!\bigcup_{n\in\omega}~\mbox{Graph}({f_n}_{\vert S_n\setminus X_n})$ is a $\Ana$ relation on $X$ with uncountable Borel chromatic number.\bigskip

 By 4D.2 and 4D.14 in [M], $D_X\!:=\!\{x\!\in\! X\mid x\!\in\! {\it\Delta }^1_1\}$ is countable and $\Ca$. We set 
$$S\!:=\!\{ x\!\in\! P\mid (x,x)\!\in\!\overline{A'}^{{\it\Delta}_X^2}\}\cap\Omega_X\!\setminus\! D_X\mbox{,}$$ 
$\tau\! :=\! {{\it\Sigma }_X}_{\vert S}$ and, for each $n\!\in\!\omega$, 
$C_n\! :=\! S\cap f_n^{-1}(S)\cap S_n\!\setminus\! X_n$.

\vfill\eject

 Note that $S$ is a Borel and $\Ana$ subset of $X$. Thus the $C_n$'s are $\Ana$ and clopen subsets of $Y$. As $S\!\in\!\Ana$, $Y$ is a zero-dimensional perfect Polish space. We set $C\! :=\!\{ x\!\in\! P\mid (x,x)\!\in\!\overline{A'}^{{\it\Delta}_X^2}\}$.\bigskip

 Let us check that $C\!\setminus\! D_X$ is nonempty. We argue by contradiction. Then $C$ is $\Ana$ and contained in $D_X$. The effective separation result gives $\Delta\!\in\!\Borel (X)$ with $C\!\subseteq\!\Delta\!\subseteq\! D_X$. If $x\!\in\! P\!\setminus\!\Delta$, then there is $U_x\!\in\!\Borel (X)$ containing $x$ with $A'\cap U_x^2\! =\!\emptyset$. This also holds if $x\!\in\!\Delta$ with $U_x\! :=\!\{ x\}$. But this contradicts the fact that $A'$ has uncountable Borel chromatic number.\bigskip

 Thus $C\!\setminus\! D_X$ is a nonempty $\Ana$ subset of $X$, which therefore meets $\Omega_X$. This shows that $S$ is not empty. Note also that 
$(x,x)\!\in\!\overline{A'}^{{\it\Delta}_X^2}\cap S^2\! =\!\overline{A'}^{{\it\Sigma}_X^2}\cap S^2\! =\!
\overline{A'\cap S^2}^{Y^2}$ if $x\!\in\! S$. We proved that $\big( Y,({f_n}_{\vert C_n})\big)$ is a strongly complex situation.\hfill{$\square$}

\begin{cor} \label{charmin+} Let $\big(X,(f_{n})\big)$ be a complex situation such that 
$(2^\omega ,\mathbb{G}_0)\not\preceq_c^{\mbox{inj}}(X,A^f)$. The following are equivalent:\smallskip

(a) $A^f$ is $\preceq^{\mbox{inj}}_c$-minimal among analytic digraphs of uncountable Borel chromatic number,\smallskip
  
(b) for any Borel subset $S$ of $X$, any finer topology $\tau$ on $S$, and any sequence $(C_n)$ of clopen subsets of $Y\! :=\! (S,\tau )$, 
$(X,A^f)\preceq^{\mbox{inj}}_c\big( Y,\bigcup_{n\in\omega}~\mbox{Graph}({f_n}_{\vert C_n})\big)$ if $\big( Y,({f_n}_{\vert C_n})\big)$ is a strongly complex situation.\end{cor} 

\noindent {\bf Proof.} (a) $\Rightarrow$ (b) By Corollary \ref{D2+}, 
$\bigcup_{n\in\omega}~\mbox{Graph}({f_n}_{\vert C_n})$ is an analytic digraph of uncountable Borel chromatic number. As $( Y,\bigcup_{n\in\omega}~\mbox{Graph}({f_n}_{\vert C_n})\big)\preceq^{\mbox{inj}}_c(X,A^f)$, we are done.\bigskip

\noindent (b) $\Rightarrow$ (a) Let $Z$ be a Polish space, and $A$ be an analytic digraph on $Z$ of uncountable Borel chromatic number. We assume that 
$(Z,A)\preceq^{\mbox{inj}}_B(X,A^f)$, with witness $u\! :\! Z\!\rightarrow\! X$. We set $P\! :=\! u[Z]$, so that $P$ is a Borel subset of $X$. Note that 
$A'\! :=\! (u\!\times\! u)[A]\!\subseteq\! A^f$ is an analytic relation on $P$ with uncountable Borel chromatic number, which gives a sequence 
$(S_n)_{n\in\omega}$ of analytic subsets of $X$ with $A'\! =\!\bigcup_{n\in\omega}~\mbox{Graph}({f_n}_{\vert S_n})$. Lemma \ref{Bas0+} gives a Borel subset $T$ of $P$, a finer topology $\sigma$ on $T$, and a sequence $(O_n)$ of clopen subsets of $W\! :=\! (T,\sigma )$ such that 
$\big( W,({f_n}_{\vert O_n})\big)$ is a strongly complex situation and $O_n\!\subseteq\! T\cap f_n^{-1}(T)\cap S_n$. We put 
$A''\! :=\!\bigcup_{n\in\omega}~\mbox{Graph}({f_n}_{\vert O_n})$.\bigskip

 Note that $b\! :=\! {u^{-1}}_{\vert W}\! :\! W\!\rightarrow\! u^{-1}(W)$ is Borel and one-to-one. Let $G$ be a dense $G_{\delta}$ subset of $W$ such that $b\vert_G$ is continuous. This function is a witness for the fact that $(G,A''\cap G^2)\preceq^{\mbox{inj}}_c(Z,A)$. By Lemma \ref{suffnpc}, we get 
$\mbox{Graph}({f_n}_{\vert O_n})\!\subseteq\!\overline{\mbox{Graph}({f_n}_{\vert O_n})\cap G^2}$. Thus 
$\overline{A''}\! =\!\overline{A''\cap G^2}$, and $\Delta (G)\!\subseteq\!\overline{A''\cap G^2}$. Let us prove that $A''\cap G^2$ has uncountable Borel chromatic number. We argue by contradiction, which gives a countable partition $(B_q)$ of $G$ into Borel sets. We can find $q\!\in\!\omega$, a nonempty open subset $O$ of 
$W$, and a dense $G_\delta$ subset $H$ of $W$ with $O\cap H\!\subseteq\! B_q$. The previous argument shows that $\Delta (G\cap H)\!\subseteq\!\overline{A''\cap (G\cap H)^2}$. Therefore 
$\Delta (O\cap G\cap H)\!\subseteq\!\overline{A''\cap (O\cap G\cap H)^2}\!\subseteq\!
\overline{A''\cap G^2\cap B_q^2}\! =\!\emptyset$, which is absurd.\bigskip

 We apply Lemma \ref{Bas0+} to $\big( W,({f_n}_{\vert O_n})\big)$, $G$ and $\big( G\cap f_n^{-1}(G)\cap O_n\big)$, which gives a Borel subset $S$ of $G$, a topology $\tau$ on $S$ finer than $\sigma$, and a sequence $(C_n)$ of clopen subsets of $Y\! :=\! (S,\tau )$ such that $\big( Y,({f_n}_{\vert C_n})\big)$ is a strongly complex situation and $C_n\!\subseteq\! S\cap f_n^{-1}(S)\cap O_n$. It remains to note that 
${(X,A^f)\!\preceq^{\mbox{inj}}_c\!\big( Y,\bigcup_{n\in\omega}\mbox{Graph}({f_n}_{\vert C_n})\big)\!\preceq^{\mbox{inj}}_c\! (G,A''\cap G^2)\preceq^{\mbox{inj}}_c(Z,A)}$, by (b).
\hfill{$\square$}

\vfill\eject

\noindent {\bf Remark.} This proof also shows that (b) implies that $A^f$ is $\preceq^{\mbox{inj}}_B$-minimal among analytic digraphs of uncountable Borel chromatic number.\bigskip

\noindent {\bf Proof of Theorem \ref{fav}.} By Corollary \ref{D2+}, (a) and (b) cannot hold simultaneously. So assume that (a) does not hold. Theorem \ref{G0} provides $h\! :\! 2^\omega\!\rightarrow\! Y$ continuous with 
$(h\!\times\! h)[\mathbb{G}_0]\!\subseteq\! B$. Note that $\mathbb{G}_0$ comes from a complex situation. By Corollary \ref{D2+}, $\mathbb{G}_0$ is a $\boratwo$ relation on the compact space $2^\omega$. Thus 
$\mathbb{G}_0$ and $(h\!\times\! h)[\mathbb{G}_0]$ are $K_\sigma$ digraphs of uncountable Borel chromatic number. The canonical injection is a witness for the fact that 
$\big( Y,(h\!\times\! h)[\mathbb{G}_0]\big)\preceq_c^{\mbox{inj}}(Y,B)$. So we may assume that $B$ is $K_\sigma$. Let $(C_n)$ be a sequence of closed digraphs on $Y$ whose union is $B$.\bigskip

  We now argue essentially as in the proof of Theorem 2.3 in [L2]. In order to simplify the notation, we will assume as usual that $Y$ is recursively presented and $(C_n)$ is $\Borel$. As in the proof of Theorem \ref{G0}, we set $X_1\! :=\! Y\!\setminus\!\big(\bigcup\{ S\!\in\!\Ana (Y)\mid B\cap S^2\! =\!\emptyset\}\big)$, so that $X_1$ is a nonempty $\Ana$ subset of $Y$, disjoint from $\{ y\!\in\! Y\mid y\!\in\! {\it\Delta }^1_1\}$, and satisfies the following property:
$$\forall S\!\in\!\Ana (Y)~~(\emptyset\!\not=\! S\!\subseteq\! X_1\Rightarrow B\cap S^2\!\not=\!\emptyset ).$$ 
We set $Z\! :=\! (X_1\cap\Omega_Y,{\it\Sigma}_Y)$, so that $Z$ is a nonempty zero-dimensional perfect Polish space. We also set ${\mathfrak S}\! :=\!\{ S\!\in\!\Ana (Z)\mid\emptyset\!\not=\! S\!\subseteq\! Z\}$. If $S\!\in\!\mathfrak{S}$, then we can find 
$n$ with $C_n\cap S^2\!\not=\!\emptyset$. Note that $C_n\cap S^2$ is a closed relation on $Z$. Moreover, if $U$ is an open relation on $Z$, then the projections of $U\cap C_n\cap S^2$ are open subsets of $Z$.\bigskip

 Theorem 1.13 in [L2] provides, for each $S\!\in\!\mathfrak{S}$, dense $G_\delta$ subsets $F_S$, $G_S$ of some nonempty open subsets of $Z$ and $g_S\! :\! F_S\!\rightarrow\! G_S$ onto, continuous and open such that 
$\mbox{Graph}(g_S)\!\subseteq\! C_n\cap S^2$ or $\mbox{Graph}(g_S)\!\subseteq\! (C_n\cap S^2)^{-1}$. We set 
$$G(g_S)\! :=\!\left\{\!\!\!\!\!\!\!\!
\begin{array}{ll} 
& \mbox{Graph}(g_S)\mbox{ if }\mbox{Graph}(g_S)\!\subseteq\! C_n\cap S^2\mbox{,}\cr
& \big(\mbox{Graph}(g_S)\big)^{-1}\mbox{ if }\mbox{Graph}(g_S)\!\subseteq\! (C_n\cap S^2)^{-1}\mbox{,}
\end{array}
\right.$$
so that $G(g_S)\!\subseteq\! B\cap S^2$.\bigskip

 Let $\alpha\!\in\! 2^\omega$ such that $Z$ is recursively in $\alpha$ presented and the sequence 
$\big( G(g_S)\big)_{S\in\mathfrak{S}}$ is $\Borel (\alpha )$. Let $G$ be a dense $G_\delta$ and $\Ana (\alpha )$ subset of 
$Z$ on which ${\it\Sigma}_Y$ and ${\it\Delta}^\alpha_Z\! :=<\Borel (\alpha )(Z)>$ coincide, which exists by Lemma 2.1 in [L2]. Let ${\it\Sigma}^\alpha_Z\! :=<\Ana (\alpha )(Z)>$, which gives $\Omega^\alpha_Z$ like before Lemma \ref{Bas0+}. Moreover, the proof of Theorem 2.3 in [L2] shows that $\Omega^\alpha_Z$ is comeager in $Z$, as well as $G\cap\Omega^\alpha_Z$. We set $W\! :=\! (G\cap\Omega^\alpha_Z,{\it\Sigma}^\alpha_Z)$ and define, for $S\!\in\!\mathfrak{S}$, a partial function $f_S$ by $G(f_S)\! :=\! G(g_S)\cap W^2$, so that $W$ is a nonempty zero-dimensional Polish space and 
$\big( W,\bigcup_{S\in\mathfrak{S}}~G(f_S)\big)\preceq_c^{\mbox{inj}}(Y,B)$. Moreover, for each $S\!\in\!\mathfrak{S}$, $f_S$ is a partial continuous and open map with clopen domain and range in $W$. We then note that 
$$\Delta (W)\!\subseteq\!\overline{\bigcup_{S\in\mathfrak{S}}~G(g_S)}^{Z^2}\cap W^2\! =\!
\overline{\bigcup_{S\in\mathfrak{S}}~G(g_S)}^{({\it\Delta}^\alpha_Z)^2}\cap W^2\! =\!
\overline{\bigcup_{S\in\mathfrak{S}}~G(g_S)}^{({\it\Sigma}^\alpha_Z)^2}\cap W^2\!\subseteq\!
\overline{\bigcup_{S\in\mathfrak{S}}~G(f_S)}^{W^2}.$$ 
We set $\mathfrak{S}_0\! :=\!\{ S\!\in\!\mathfrak{S}\mid\mbox{Graph}(g_S)\!\subseteq\! B\}$, 
$\mathfrak{S}_1\! :=\!\{ S\!\in\!\mathfrak{S}\mid\mbox{Graph}(g_S)\!\subseteq\! B^{-1}\}$ and, for $\varepsilon\!\in\! 2$, 
$$W_\varepsilon\! :=\!\{ w\!\in\! W\mid (w,w)\!\in\!\overline{\bigcup_{S\in\mathfrak{S}_\varepsilon}~G(f_S)}\} .$$ 

 Note that $(W_0,W_1)$ is a covering of $W$ into closed sets, which gives a nonempty clopen subset $X$ of $W$ and $\varepsilon\!\in\! 2$ with $X\!\subseteq\! W_\varepsilon$. Note that $X$ is a nonempty Polish space, $G(f_S)\cap X^2$ defines a partial continuous and open map with open domain and range in $X$ if $S\!\in\!\mathfrak{S}_\varepsilon$, and 
$$\Delta (X)\!\subseteq\!\overline{\bigcup_{S\in\mathfrak{S}_\varepsilon}~G(f_S)\cap X^2}^{X^2}\!\setminus\!
\big(\bigcup_{S\in\mathfrak{S}_\varepsilon}~G(f_S)\cap X^2\big) .$$ 
We enumerate $(f_n)\! :=\!\big( (f_S)_{\vert X\cap f_S^{-1}(X)}\big)_{S\in\mathfrak{S}_\varepsilon}$. If $\varepsilon\! =\! 0$, then $\big( X,(f_n)\big)$ is a complex situation with $(X,A^f)\preceq_c^{\mbox{inj}}(Y,B)$. If $\varepsilon\! =\! 1$, then 
$\big( X,(f_n)\big)$ is a complex situation with $(X,A^f)\preceq_c^{\mbox{inj}}(Y,B^{-1})$.\hfill{$\square$}

\section{$\!\!\!\!\!\!$ The lemma ensuring the injectivity}\indent

 Recall the set $X_l$ mentioned at the end of the introduction. An oriented graph $A_l$ on $X_l$ will contain finite approximations of a subset of $\mathbb{G}_1$ with acyclic symmetrization. We now isolate some important properties of $A_l$ leading to a relatively general lemma ensuring the injectivity of the homomorphism in Theorem \ref{cormain}.\bigskip

\noindent\bf Notation.\rm\ Let $X$ be a set, and $A$ be an oriented graph on $X$. We set, for $x\!\in\! X$, 
$$\left\{\!\!\!\!\!\!\!
\begin{array}{ll}
& \mbox{Succ}(x)\! :=\!\{ y\!\in\! X\mid (x,y)\!\in\! A\}\mbox{ and }
\mbox{max}_X\! :=\!\{ x\!\in\! X\mid\mbox{Succ}(x)\! =\!\emptyset\}\mbox{,}\cr\cr
& \mbox{Pred}(x)\! :=\!\{ y\!\in\! X\mid (y,x)\!\in\! A\}\mbox{ and }
\mbox{min}_X\! :=\!\{ x\!\in\! X\mid\mbox{Pred}(x)\! =\!\emptyset\} .
\end{array}
\right.$$
The following oriented graphs will be of particular importance in the sequel.

\begin{defi} An oriented graph $A$ on a set $X$ is {\bf unambiguously oriented} if 
$\vert\mbox{Succ}(x)\vert\!\leq\! 1$ for each $x\!\in\! X$. If moreover $A$ has acyclic symmetrization, then we say that $A$ is an {\bf uogas}.\end{defi}

\begin{lem} \label{am1} Let $A$ be an uogas on a finite set $X$.\smallskip

(a) Let $l$ be a natural number, $(y_i)_{i\leq l}\!\in\! X^{l+1}$ such that $y_{i+1}\!\in\!\mbox{Succ}(y_i)$ if 
$i\! <\! l$. Then $(y_i)_{i\leq l}$ is injective. In particular, $(y_i)_{i\leq l}\! =\! p^{s(A)}_{y_0,y_l}$.\smallskip

(b) Let $y\!\in\! X$, $x\!\in\! C(y)\cap\mbox{max}_X$, and $p\! :=\! p^{s(A)}_{y,x}$. Then 
$\big( p(i),p(i\! +\! 1)\big)\!\in\! A$ if $i\! <\!\vert p\vert\! -\! 1$. \smallskip

(c) The intersection of $\mbox{max}_X$ with each $s(A)$-connected component $C$ is a singleton 
$\{ x_C\}$.\smallskip

(d) Let $(y,x)\!\in\! A$ and $p\! :=\! p^{s(A)}_{y,x_{C(y)}}$. Then $\vert p\vert\!\geq\! 2$ and 
$p(1)\! =\! x$.\end{lem}

\noindent\bf Proof.\rm\ (a) It is enough to see that $(y_i)_{i\leq l}$ is injective. This is clear if $l\! =\! 0$. As $A$ is irreflexive, $y_i\!\not=\! y_{i+1}$ if $i\! <\! l$. As $A$ is an oriented graph, 
$y_i\!\not=\! y_{i+2}$ if $i\! <\! l\! -\! 1$. As $s(A)$ is acyclic, $(y_i)_{i\leq l}$ is injective, by induction on $l$.\bigskip

\noindent (b) We argue by induction on $m\! :=\!\vert p\vert\! -\! i$. For $m\! =\! 2$, we argue by contradiction, so that $\big( x,p(i)\big)\!\in\! A$. Thus $p(i)\!\in\!\mbox{Succ}(x)$, which contradicts the fact that $x\!\in\!\mbox{max}_X$. So assume that 
$$\big( p(\vert p\vert\! -\! m),p(\vert p\vert\! -\! m\! +\! 1)\big)\!\in\! A.$$ 
We argue by contradiction, so that $\big( p(\vert p\vert\! -\! m),p(\vert p\vert\! -\! m\! -\! 1)\big)\!\in\! A$. Thus $p(\vert p\vert\! -\! m\! +\! 1)$ and $p(\vert p\vert\! -\! m\! -\! 1)$ are different elements of 
$\mbox{Succ}\big( p(\vert p\vert\! -\! m)\big)$, which contradicts the fact that $A$ is unambiguously oriented.

\vfill\eject

\noindent (c) Let $C$ be a $s(A)$-connected component, and $\{ x_m\mid m\! <\! l\}$ be an injective enumeration of $C$. We argue by contradiction to prove that $\mbox{max}_X$ meets $C$. This means that 
$\mbox{Succ}(x_m)\!\not=\!\emptyset$ if $m\! <\! l$. We define, inductively on $i$, a sequence 
$(m_i)_{i\leq l}$ as follows. We first set $m_0\! :=\! 0$. Assume that $i\! <\! l$ and $m_i$ has been defined. As $\mbox{Succ}(x_{m_i})\!\not=\!\emptyset$, there is $m_{i+1}\! <\! l$ such that 
$x_{m_{i+1}}\!\in\!\mbox{Succ}(x_{m_i})$. By (a), 
$(m_i)_{i\leq l}\!\in\! l^{l+1}$ is injective, which is absurd.\bigskip

 Assume now that, for example, $x_0,x_1\!\in\! C\cap\mbox{max}_X$. Let $p_0\! :=\! p^{s(A)}_{x_0,x_1}$ and $p_1\! :=\! p^{s(A)}_{x_1,x_0}$, so that $p_1(i)\! =\! p_0(n\! -\! i\! -\! 1)$ if 
$i\! <\! n\! :=\!\vert p_1\vert\! =\!\vert p_0\vert$. By (b), 
$\big( p_0(n\! -\! 2),x_1\big) ,\big( x_1,p_1(1)\big)\!\in\! A$, which contradicts the fact that $A$ is an oriented graph.\bigskip

\noindent (d) We set $y_0\! :=\! y$, $y_1\! :=\! x$, and choose $y_{i+1}\!\in\!\mbox{Succ}(y_i)$ if this last set is not empty. After finitely many steps, this construction stops and provides $(y_i)_{i\leq l}$ with $l\!\geq\! 1$, by (a) and since $X$ is finite. Note that $y_l\! =\! x_{C(y)}$ and $(y_i)_{i\leq l}\! =\! p$, by construction and by (a) and (c).\hfill{$\square$}\bigskip

 The homomorphism $h$ in Theorem \ref{cormain} will be obtained thanks to a Cantor-like construction. In the inductive step of this construction, we will consider an uogas on a finite set $X$. We will associate open sets to the elements of $X$. As we want $h$ to be injective, we will have to ensure the disjunction of these open sets. This will be achieved in Lemma \ref{disj} to come. Its proof will use the following objects.\bigskip

\noindent\bf Notation.\rm\ Let $A$ be an uogas on a finite set $X$. If $y\!\in\! X$, then Lemma \ref{am1} allows us to set $p_y\! :=\! p^{s(A)}_{y,x_{C(y)}}$. We set $L\! :=\!\vert X\vert$, and enumerate 
${X\! :=\!\{ x_m\vert m\! <\! L\}}$ injectively in such a way that $(\vert p_{x_m}\vert )_{m<L}$ is increasing. Let $L_0\!\leq\! L$ be maximal such that $x_m\!\in\!\mbox{max}_X$ if $m\! <\! L_0$, so that 
$\{ x_m\vert m\! <\! L_0\}$ enumerates injectively $\mbox{max}_X$. The idea is to make, inductively on 
$m$, $L$ copies of $x_m$ if $m\!\geq\! L_0$, keeping the $A$-relations. In order to do this, we will give labels to the elements of $X$. We set 
$${\cal X}_m\! :=\!\big\{\big( x_k,(0)\big)\mid k\! <\! L\big\}$$ 
and ${\cal A}_m\! :=\!\Big\{\Big(\big( x_k,(0)\big) ,\big( x_l,(0)\big)\Big)\mid (x_k,x_l)\!\in\! A\Big\}$ if 
$m\! <\! L_0$. Assume that $L_0\!\leq\! m\! +\! 1\! <\! L$. We put the elements of ${\cal X}_m$ whose $X$-coordinate is not $x_{m+1}$ or one of its (iterated) predecessors in 
$${\cal R}_{m+1}\! :=\!\big\{ (x_k,\sigma )\!\in\! {\cal X}_m\mid x_{m+1}\!\notin\! p_{x_k}\big\} .$$ 
We also set, for $\sigma\!\in\! L^{\leq m+1}$ and $j\! <\! L$, ${\cal X}^\sigma_{m+1}\! :=\!
\big\{ (x_k,\sigma )\!\in\! {\cal X}_m\mid x_{m+1}\!\in\! p_{x_k}\big\}$, 
$${\cal X}^{\sigma j}_{m+1}\! :=\!\big\{ (x_k,\sigma j)\!\in\! X\!\times\! L^{<\omega}\mid 
(x_k,\sigma )\!\in\! {\cal X}_m\wedge x_{m+1}\!\in\! p_{x_k}\big\}$$  
(some of these sets can be empty). Then we set 
${\cal X}_{m+1}\! :=\! {\cal R}_{m+1}\cup\bigcup_{\sigma\in L^{\leq m+1},j<L}~{\cal X}^{\sigma j}_{m+1}$ and\bigskip

\leftline{${\cal A}_{m+1}\! :=\!\big\{\big( (x_k,\sigma ) ,(x_l,\tau )\big)\!\in\! {\cal A}_m\mid k,l\!\leq\! m\big\}\ \cup$}\smallskip

\rightline{$\big\{\big( (x_{m+1},\sigma j) ,(x_i,\sigma )\big)\!\in\! ({\cal X}_{m+1})^2\mid 
\big( (x_{m+1},\sigma ) ,(x_i,\sigma)\big)\!\in\! {\cal A}_m\wedge j\! <\! L\big\}\ \cup$}\smallskip

\rightline{$\big\{\big( (x_k,\sigma ) ,(x_l,\sigma )\big)\!\in\! {\cal A}_m\mid k\! >\! m\! +\! 1\wedge 
x_{m+1}\!\notin\! p_{x_l}\big\}\ \cup$}\smallskip

\rightline{$\big\{\big( (x_k,\sigma j) ,(x_l,\sigma j)\big)\!\in\! ({\cal X}_{m+1})^2\mid 
\big( (x_k,\sigma ) ,(x_l,\sigma)\big)\!\in\! {\cal A}_m\wedge j\! <\! L\wedge x_{m+1}\!\in\! p_{x_l}\big\} .$}\bigskip

\noindent We enumerate injectively $\{\sigma\!\in\! L^{<\omega}\mid (x_{m+1},\sigma )\!\in\! {\cal X}_m\}$ by 
$\{\sigma_n\mid n\! <\! N\}$. We set, for $p\!\leq\! N$, 
${\cal X}^p_{m+1}\! :=\! {\cal R}_{m+1}\cup\bigcup_{n<p,j<L}~{\cal X}^{\sigma_nj}_{m+1}\cup
\bigcup_{p\leq n<N}~{\cal X}^{\sigma_n}_{m+1}$, so that ${\cal X}^0_{m+1}\! =\! {\cal X}_m$ and 
${\cal X}^N_{m+1}\! =\! {\cal X}_{m+1}$. We also define the corresponding intermediate versions of ${\cal A}_{m+1}$ as follows.

\vfill\eject

 We set, for $p\!\leq\! N$,\bigskip

\leftline{${\cal A}^p_{m+1}\! :=\!\big\{\big( (x_k,\sigma ) ,(x_l,\tau )\big)\!\in\! {\cal A}_m\mid k,l\!\leq\! m\big\}\ \cup$}\smallskip

\rightline{$\bigcup_{n<p}~\big\{\big( (x_{m+1},\sigma_nj) ,(x_i,\sigma_n)\big)\!\in\! ({\cal X}^p_{m+1})^2\mid 
\big( (x_{m+1},\sigma_n) ,(x_i,\sigma_n)\big)\!\in\! {\cal A}_m\wedge j\! <\! L\big\}\ \cup$}\smallskip

\rightline{$\bigcup_{p\leq n<N}~\big\{\big( (x_{m+1},\sigma_n) ,(x_i,\sigma_n)\big)\!\in\! 
({\cal X}^p_{m+1})^2\mid\big( (x_{m+1},\sigma_n) ,(x_i,\sigma_n)\big)\!\in\! {\cal A}_m\big\}\ \cup$}\smallskip

\rightline{$\big\{\big( (x_k,\sigma ) ,(x_l,\sigma )\big)\!\in\! {\cal A}_m\mid 
k\! >\! m\! +\! 1\wedge x_{m+1}\!\notin\! p_{x_l}\big\}\ \cup$}\smallskip

\rightline{$\bigcup_{n<p}~\big\{\big( (x_k,\sigma_nj) ,(x_l,\sigma_nj)\big)\!\in\! ({\cal X}^p_{m+1})^2\mid 
\big( (x_k,\sigma_n) ,(x_l,\sigma_n)\big)\!\in\! {\cal A}_m\wedge j\! <\! L\wedge x_{m+1}\!\in\! p_{x_l}\big\}\ \cup$}\smallskip

\rightline{$\bigcup_{p\leq n<N}~\big\{\big( (x_k,\sigma_n) ,(x_l,\sigma_n)\big)\!\in\! {\cal A}_m\mid 
x_{m+1}\!\in\! p_{x_l}\big\} .$}

\begin{lem} \label{split} Let $A$ be an uogas on a finite set $X$. Then ${\cal A}_m$ (respectively, 
${\cal A}^p_{m+1}$) is also an uogas on the finite set ${\cal X}_m$ (respectively, ${\cal X}^p_{m+1}$) if $m\! <\! L$ (respectively, $L_0\!\leq\! m\! +\! 1\! <\! L$ and $p\!\leq\! N$).\end{lem}

\noindent\bf Proof.\rm\ Note first that $({\cal X}_m,{\cal A}_m)$ is a copy of $(X,A)$ if $m\! <\! L_0$. In particular, ${\cal X}_m$ is finite, ${\cal A}_m$ is an oriented graph on ${\cal X}_m$ with acyclic symmetrization, and $\vert\mbox{Succ}(q)\vert\!\leq\! 1$ for each $q\!\in\! {\cal X}_m$ (with respect to ${\cal A}_m$). Then we assume that $L_0\!\leq\! m\! +\! 1\! <\! L$. Note that ${\cal X}_{m+1}$ is finite, ${\cal A}_{m+1}$ is an oriented graph on ${\cal X}_{m+1}$, and 
$\vert\mbox{Succ}(q)\vert\!\leq\! 1$ for each $q\!\in\! {\cal X}_{m+1}$ (with respect to ${\cal A}_{m+1}$). Let us check that $s({\cal A}_{m+1})$ is acyclic. The restriction of ${\cal A}_{m+1}$ to a fixed ${\cal X}^{\sigma j}_{m+1}$ is isomorphic to a subgraph of ${\cal A}_m$, and has therefore acyclic symmetrization. Note that the ${\cal X}^{\sigma j}_{m+1}$'s are pairwise disjoint and not $s({\cal A}_{m+1})$-related, and that ${\cal X}^{\sigma j}_{m+1}$ contains 
$(x_{m+1},\sigma j)$ if it is not empty (i.e., if $(x_{m+1},\sigma )$ is in ${\cal X}_m$). The restriction of the oriented graph ${\cal A}_{m+1}$ to ${\cal R}_{m+1}$ is also isomorphic to a subgraph of ${\cal A}_m$. Moreover, the only possible $s({\cal A}_{m+1})$-edge between an element of ${\cal X}^{\sigma j}_{m+1}$ and element of ${\cal R}_{m+1}$ is between $(x_{m+1},\sigma j)$ and $(x_i,\sigma )$, where $x_i\!\in\!\mbox{Succ}(x_{m+1})$. This shows the acyclicity of 
$s({\cal A}_{m+1})$. We argue similarly for $({\cal X}^p_{m+1},{\cal A}^p_{m+1})$.\hfill{$\square$}

\begin{defi} We say that the tuple ${\cal T}\! :=\! \big( X,A,Z,(f_n)\big)$ is a {\bf mapping tuple} if $A$ is an uogas on the finite set $X$ and $\big( Z,(f_n)\big)$ is a strongly complex situation.\end{defi}

\noindent\bf Notation.\rm\ Let ${\cal T}\! :=\!\big( X,A,Z,(f_n)\big)$ be a mapping tuple. We set 
$$\left\{\!\!\!\!\!\!\!
\begin{array}{ll}
& E_{\cal T}\! :=\!
\big\{ ~\big( u,(V_x)_{x\in X}\big)\!\in\!\omega^X\!\times\! (\boraone (Z)\!\setminus\!\{\emptyset\} )^X\mid
\forall (x,y)\!\in\! A~~V_x\!\subseteq\! D_{u(x)}\wedge V_y\! =\! f_{u(x)}[V_x]~\big\}\mbox{,}\cr\cr
& U_{\cal T}\! :=\!
\big\{ ~\big( u,(V_x)_{x\in X}\big)\!\in\!\omega^X\!\times\! (\boraone (Z)\!\setminus\!\{\emptyset\} )^X\mid
\forall (x,y)\!\in\! A~~V_x\!\subseteq\! D_{u(x)}\wedge V_y\!\subseteq\! f_{u(x)}[V_x]~\big\} .
\end{array}
\right.$$

\begin{lem} \label{equal} Let ${\cal T}\! :=\!\big( X,A,Z,(f_n)\big)$ be a mapping tuple, and 
$\big( u,(V_x)_{x\in X}\big)\!\in\! U_{\cal T}$. Then we can find a family $(W_x)_{x\in X}$ of subsets of $Z$ such that\smallskip

(a) $\big( u,(W_x)_{x\in X}\big)\!\in\! E_{\cal T}$,\smallskip

(b) $W_x\!\subseteq\! V_x$ if $x\!\in\! X$, and $W_x\! =\! V_x$ if $x\!\in\!\mbox{max}_X$.\end{lem}

\noindent\bf Proof.\rm\ We define $W_x$ by induction on $\vert p_x\vert$. If $\vert p_x\vert\! =\! 1$, then we set $W_x\! :=\! V_x$. Assume that 
$\vert p_x\vert\!\geq\! 2$, so that $W_y$ has been defined if $y\! :=\! p_x(1)$. We set $W_x\! :=\! V_x\cap f_{u(x)}^{-1}(W_y)$. We are done, by Lemma \ref{am1}.\hfill{$\square$}

\begin{lem} \label{propag} Let ${\cal T}\! :=\!\big( X,A,Z,(f_n)\big)$ be a mapping tuple, 
$\big( u,(V_x)_{x\in X}\big)\!\in\! E_{\cal T}$, $x_0\!\in\! X$, and $W_{x_0}$ be a nonempty open subset of $V_{x_0}$. Then we can find a family $(W_x)_{x\in X\setminus\{ x_0\}}$ of subsets of $Z$ such that\smallskip

(a) $\big( u,(W_x)_{x\in X}\big)\!\in\! E_{\cal T}$,\smallskip

(b) $W_x\!\subseteq\! V_x$ if $x\!\in\! X$, and $W_x\! =\! V_x$ if $x\!\notin\! C(x_0)$.\end{lem}

\vfill\eject

\noindent\bf Proof.\rm\ We set, for $x\!\in\! C(x_0)$, $q_x\! :=\! p^{s(A)}_{x,x_0}$. We define $W_x$ by induction on $\vert q_x\vert$, the case 
$\vert q_x\vert\! =\! 1$ (i.e., $x\! =\! x_0$) being done. So assume that $\vert q_x\vert\!\geq\! 2$, so that $W_y$ has been defined if $y\! :=\! q_x(1)$. We set 
$$W_x\! :=\!\left\{\!\!\!\!\!\!\!\!
\begin{array}{ll}
& f_{u(y)}[W_y]\mbox{ if }(y,x)\!\in\! A\mbox{,}\cr\cr
& V_x\cap f_{u(x)}^{-1}(W_y)\mbox{ if }(x,y)\!\in\! A.
\end{array}
\right.$$
If $x\!\in\! X\!\setminus\! C(x_0)$, then we set $W_x\! :=\! V_x$.\hfill{$\square$}\bigskip

\noindent\bf Notation.\rm\ Let $A$ be an uogas on $X$ finite, and $x\!\in\! X$. We set 
$M_x\! :=\!\mbox{max}\{\vert p^{s(A)}_{y,x}\vert\mid y\!\in\!\mbox{min}_X\wedge x\!\in\! p_y\}$. 

\begin{lem} \label{disj} Let ${\cal T}\! :=\!\big( X,A,Z,(f_n)\big)$ be a mapping tuple, $d\!\in\!\omega$, and 
$\big( u,(V_x)_{x\in X}\big)\!\in\! U_{\cal T}$. Then we can find a family $(W_x)_{x\in X}$ of subsets of $Z$ such that\smallskip

(a) $\big( u,(W_x)_{x\in X}\big)\!\in\! E_{\cal T}$,\smallskip

(b) $W_x\!\subseteq\! V_x$, $W_x\!\in\!\borone (Z)$ and $\mbox{diam}(W_x)\!\leq\! 2^{-d}$ if $x\!\in\! X$,\smallskip

(c) $W_x\cap W_y\! =\!\emptyset$ if $x\!\not=\! y\!\in\! X$.\end{lem}

\noindent\bf Proof.\rm\ We consider the oriented graphs ${\cal A}_m$ (respectively, ${\cal A}^p_{m+1}$) on 
${\cal X}_m$ (respectively, ${\cal X}^p_{m+1}$) defined before Lemma \ref{split}. We set, for $m\! <\! L$, 
${\cal T}_m\! :=\!\big( {\cal X}_m,{\cal A}_m,Z,(f_n)\big)$. Similarly, we set, when ${L_0\!\leq\! m\! +\! 1\! <\! L}$ and $p\!\leq\! N$, 
${\cal T}^p_{m+1}\! :=\!\big( {\cal X}^p_{m+1},{\cal A}^p_{m+1},Z,(f_n)\big)$. By Lemma \ref{split}, all these tuples are mapping tuples. We define, for 
$m\! <\! L$, $u_m\! :\! {\cal X}_m\!\rightarrow\!\omega$ by $u_m(x_k,\sigma )\! :=\! u(x_k)$. We also define, when $L_0\!\leq\! m\! +\! 1\! <\! L$ and $p\!\leq\! N$, $u^p_{m+1}\! :\! {\cal X}^p_{m+1}\!\rightarrow\!\omega$ by $u^p_{m+1}(x_k,\sigma )\! :=\! u(x_k)$.\bigskip

 We first construct, by induction on $m\! <\! L$, a family $(W^m_q)_{q\in {\cal X}_m}$ of nonempty open subsets of $Z$ satisfying
$$\begin{array}{ll} 
& (1)~W^{m+1}_{x_k,\sigma'}\!\subseteq\! W^m_{x_k,\sigma}\!\subseteq\! V_{x_k}\cr    
& (2)~\big( u_m,(W^m_q)_{q\in {\cal X}_m}\big)\!\in\! E_{{\cal T}_m}\cr
& (3)~\mbox{diam}(W^m_{x_m,\sigma})\!\leq\! 2^{-d}\cr    
& (4)~\forall j\!\not=\! j'\! <\! L~~W^m_{x_m,\sigma j}\cap W^m_{x_m,\sigma j'}\! =\!\emptyset
\end{array}$$
We will apply Lemmas \ref{equal} and \ref{propag} to perform this construction. For $m\! =\! 0$, we will apply Lemma \ref{equal} to ${\cal T}_0$ and $u_0$. We choose a nonempty open subset $V_{x_0,(0)}$ of $V_{x_0}$ with diameter at most $2^{-d}$, and we set $V_{x_k,(0)}\! :=\! V_{x_k}$ if $0\! <\! k\! <\! L$. Note that $\big( u_0,(V_q)_{q\in {\cal X}_0}\big)\!\in\! U_{{\cal T}_0}$. Lemma \ref{equal} provides a family 
$(W^0_q)_{q\in {\cal X}_0}$ of subsets of $Z$ such that\bigskip

(a) $\big( u_0,(W^0_q)_{q\in {\cal X}_0}\big)\!\in\! E_{{\cal T}_0}$,\smallskip

(b) $W^0_{x_k,(0)}\!\subseteq\! V_{x_k,(0)}\!\subseteq\! V_{x_k}$.\bigskip

\noindent This completes the construction for $m\! =\! 0$. If $m\! +\! 1\! <\! L_0$, then we proceed similarly: we choose a nonempty open subset $V_{x_{m+1},(0)}$ of $V_{x_{m+1}}$ (which is equal to $W^m_{x_{m+1},(0)}$ since we only applied Lemma \ref{equal} to perform this construction up to this point) with diameter at most $2^{-d}$, and we set $V_{x_k,(0)}\! :=\! W^m_{x_k,(0)}$ if $m\! +\! 1\!\not=\! k\! <\! L$. Note that 
$\big( u_{m+1},(V_q)_{q\in {\cal X}_{m+1}}\big)\!\in\! U_{{\cal T}_{m+1}}$. Lemma \ref{equal} provides a family $(W^{m+1}_q)_{q\in {\cal X}_{m+1}}$ of subsets of $Z$ such that\bigskip

(a') $\big( u_{m+1},(W^{m+1}_q)_{q\in {\cal X}_{m+1}}\big)\!\in\! E_{{\cal T}_{m+1}}$,\smallskip

(b') $W^{m+1}_q\!\subseteq\! W^m_q$.

\vfill\eject

 So we may assume that $L_0\!\leq\! m\! +\! 1\! <\! L$. We will, starting with 
$(Z^0_q)_{q\in {\cal X}_m}\! :=\! (W^m_q)_{q\in {\cal X}_m}$ and inductively on $p\!\leq\! N$, construct families 
$(Z^p_q)_{q\in {\cal X}^p_{m+1}}$ of subsets of $Z$ such that\bigskip

(a'') $\big( u^p_{m+1},(Z^p_q)_{q\in {\cal X}^p_{m+1}}\big)\!\in\! E_{{\cal T}^p_{m+1}}$,\smallskip

(b'') $Z^{p+1}_{x_k,\sigma'}\!\subseteq\! Z^p_{x_k,\sigma}$.\bigskip

 At the end we will set $W^{m+1}_q\! :=\! Z^N_q$. Assume that $p\! <\! N$ and 
$(Z^p_q)_{q\in {\cal X}^p_{m+1}}$ has been constructed, which is the case for $p\! =\! 0$. We will apply Lemma \ref{equal} to ${\cal T}^{p+1}_{m+1}$ and $u^{p+1}_{m+1}$. Let 
$u'\! :=\! u^p_{m+1}(x_{m+1},\sigma_p)$,  so that $Z^p_{x_{m+1},\sigma_p}\!\subseteq\! D_{u'}$ and 
$f_{u'}[Z^p_{x_{m+1},\sigma_p}]\! =\! Z^p_{x_i,\sigma_p}$, by (a''). As 
${f_{u'}}_{\vert Z^p_{x_{m+1},\sigma_p}}$ is not countable-to-one, we can find $z\!\in\! Z^p_{x_i,\sigma_p}$ and $(z_j)_{j<L}\!\in\! (Z^p_{x_{m+1},\sigma_p})^L$ injective such that $f_{u'}(z_j)\! =\! z$ for each 
$j\! <\! L$. We choose a sequence $(V_{x_{m+1},\sigma_pj})_{j<L}$ of pairwise disjoint clopen subsets of $Z$ with diameter at most $2^{-d}$ such that 
$z_j\!\in\! V_{x_{m+1},\sigma_pj}\!\subseteq\! Z^p_{x_{m+1},\sigma_p}$ for each $j\! <\! L$, and put 
$V_{x_i,\sigma_p}\! :=\! Z^p_{x_i,\sigma_p}\cap\bigcap_{j<L}~f_{u'}[V_{x_{m+1},\sigma_pj}]$, which is an open neighborhood of $z$. We apply Lemma \ref{propag} to ${\cal T}^p_{m+1}$, 
$\big( u^p_{m+1},(Z^p_q)_{q\in {\cal X}^p_{m+1}}\big)$, $(x_i,\sigma_p)$ and $V_{x_i,\sigma_p}$, which gives $\big( u^p_{m+1},(Y_q)_{q\in {\cal X}^p_{m+1}}\big)$.\bigskip

 We then set $V_{x_k,\tau}\! :=\! Y_{x_k,\tau}$ if $x_k\!\in\! p_{x_i}$ and $\tau\!\subseteq\!\sigma_p$. We also set $V_{x_k,\tau}\! :=\! Z^p_{x_k,\tau}$ if $x_k\!\notin\! p_{x_{m+1}}$, or ($x_k\!\in\! p_{x_i}$ and 
$\tau\!\not\subseteq\!\sigma_p$), or ($k\! =\! m\! +\! 1$ and $\tau$ is not of the form $\sigma_pj$ for some 
$j\! <\! L$). This defines $(V_q)_{q\in {\cal X}^{p+1}_{m+1}}$, and 
$\big( u^{p+1}_{m+1},(V_q)_{q\in {\cal X}^{p+1}_{m+1}}\big)\!\in\! U_{{\cal T}^{p+1}_{m+1}}$. Lemma \ref{equal} provides a family $(Z^{p+1}_q)_{q\in {\cal X}^{p+1}_{m+1}}$ of subsets of $Z$ such that\bigskip

(a''') $\big( u^{p+1}_{m+1},(Z^{p+1}_q)_{q\in {\cal X}^{p+1}_{m+1}}\big)\!\in\! E_{{\cal T}^{p+1}_{m+1}}$,\smallskip

(b''') $Z^{p+1}_q\!\subseteq\! V_q$ if $q\!\in\! {\cal X}^{p+1}_{m+1}$.\bigskip

\noindent This finishes the construction of the family $(W^{m+1}_q)_{q\in {\cal X}_{m+1}}$ as desired. It remains to get clopen sets and ensure (c). We first ensure the disjointness. In order to do this, we construct, inductively  on 
$m$, a family $(z_m)_{m<L}$ of points of $Z$, and a family $(\sigma_m)_{m<L}$ of finite sequences of elements of $L$ such that 
$$\begin{array}{ll} 
& (\alpha )~z_{m+1}\!\notin\!\{ z_l\mid l\!\leq\! m\}\cr    
& (\beta )~(x_m,\sigma_m)\!\in\! {\cal X}_{L-1}\cr
& (\gamma )~z_m\!\in\! W^{L-1}_{x_m,\sigma_m}\cr
& (\delta )~f_{u(x_{m+1})}(z_{m+1})\! =\! z_i\mbox{ if }x_i\!\in\!\mbox{Succ}(x_{m+1})
\end{array}$$
We first set $\sigma_0\! :=\! (0)$, so that $(x_0,\sigma_0)\!\in\! {\cal X}_{L-1}$, and choose 
$z_0\!\in\! W^{L-1}_{x_0,\sigma_0}$. If $m\! +\! 1\! <\! L_0$, then we set $\sigma_{m+1}\! :=\! (0)$, so that 
$(x_{m+1},\sigma_{m+1})\!\in\! {\cal X}_{L-1}$,  and choose 
$z_{m+1}\!\in\! W^{L-1}_{x_{m+1},\sigma_{m+1}}\!\setminus\!\{ z_l\mid l\!\leq\! m\}$, which is possible since $Z$ is perfect.\bigskip

 If $L_0\!\leq\! m\! +\! 1\! <\! L$, then let $i\!\leq\! m$ such that $x_i$ is the unique element of $\mbox{Succ}(x_{m+1})$, so that 
$(x_{m+1},\sigma_i j)\!\in\! {\cal X}_{L-1}$ for each $j\! <\! L$. We choose $j_{m+1}\! <\! L$ such that 
$\{ z_l\mid l\!\leq\! m\}\cap W^{L-1}_{x_{m+1},\sigma_i j_{m+1}}\! =\!\emptyset$, which is possible by (4). We then set $\sigma_{m+1}\! :=\!\sigma_ij_{m+1}$, so that $(x_{m+1},\sigma_{m+1})\!\in\! {\cal X}_{L-1}$,  and choose $z_{m+1}\!\in\! W^{L-1}_{x_{m+1},\sigma_{m+1}}$ with $f_{u(x_{m+1})}(z_{m+1})\! =\! z_i$.\bigskip

 Note that $\big( u,(W^{L-1}_{x_m,\sigma_m})_{m<L}\big)\!\in\! E_{\cal T}$. As $(z_m)_{m<L}$ is injective, we can find a family $(O_x)_{x\in X}$ of pairwise disjoint clopen subsets of $Z$ such that 
$z_m\!\in\! O_{x_m}\!\subseteq\! W^{L-1}_{x_m,\sigma_m}$.

\vfill\eject

 Recall the definition of $M_x$ just before Lemma \ref{disj}. We then define, for $y\!\in\! X$, and inductively on $M_y$, $U_y\! :=\! O_y\cap\bigcap_{x\in\mbox{Pred}(y)}~f_{u(x)}[U_x]$, so that 
$z_m\!\in\! U_{x_m}\!\subseteq\! V_{x_m}$, the $U_x$'s are pairwise disjoint, and 
$\big( u,(U_x)_{x\in X}\big)\!\in\! U_{\cal T}$.\bigskip

 Lemma \ref{equal} provides a family $(W_x)_{x\in X}$ of subsets of $Z$ such that\bigskip

(a'''') $\big( u,(W_x)_{x\in X}\big)\!\in\! E_{\cal T}$,\smallskip

(b'''') $W_x\!\subseteq\! U_x$.\bigskip

\noindent It remains to get clopen sets. This can be done if we apply the proof of Lemma \ref{propag}, inductively on $M_y$ as above, using the fact that this proof uses pre-images to go towards the elements of 
$\mbox{min}_X$.\hfill{$\square$}

\section{$\!\!\!\!\!\!$ The important properties of $\mathbb{G}_1$}\indent

 We first introduce the other elements of the countable antichain mentioned after Theorem \ref{biganti}.\bigskip

\noindent {\bf Notation.} Let $(g_n)$ be a sequence of partial functions. We set, for $s\!\in\!\omega^{<\omega}\!\setminus\!\{\emptyset\}$, 
$$g_s\! :=\! g_{s(0)}\circ\cdots\circ g_{s(\vert s\vert -1)}$$ 
and $s^*\! :=<s(1),\cdots ,s(\vert s\vert\! -\! 1)>$. If $L\! =\! 1$ and $j\!\geq\! 1$, then we set 
$\theta (j)\! :=\! 2j\! +\! 1$. Fix $L\!\geq\! 2$. We set $P_L\! :=\!\{ 2^p\cdot\! 3^l\mid p\!\in\!\omega\wedge l\! <\! L\! -\! 2\}$ and $M_L\! :=\!\{ 2^{31}\!\cdot\! 3\!\cdot\! k\mid k\!\geq\! 1\wedge k\!\notin\! P_L\}$. We define a map 
$\theta\! :\!\{ j\!\in\!\omega\mid j\!\geq\! 1\}\!\rightarrow\!\{ 3k\mid k\!\geq\! 1\}$ as follows:
$$\theta (j)\! :=\!\left\{\!\!\!\!\!\!\!
\begin{array}{ll}
& 3j\mbox{ if }j\!\notin\!\big( M_L\cup (M_L\! +\! 1)\big)\mbox{,}\cr
& 3j\! +\! 3\mbox{ if }j\!\in\! M_L\mbox{,}\cr
& 3j\! -\! 3\mbox{ if }j\!\in\! M_L\! +\! 1.
\end{array}
\right.$$
Recall the definition of $(S_n)$ after Theorem \ref{infi}. We define $\theta_n\! :\!\omega\!\rightarrow\!\omega$ by 
$$\theta_n(k)\! :=\!\left\{\!\!\!\!\!\!\!\!
\begin{array}{ll} 
& k\mbox{ if }k\!\notin\! S_n\mbox{,}\cr
& 2^{q_n}\!\cdot\!\theta (j)\mbox{ if }k\!\ =\! 2^{q_n}\!\cdot\! j.
\end{array}
\right.$$

\begin{lem} \label{deff++} Fix $L\!\geq\! 1$. Then $(S_n)$ and $(\theta_n)$ satisfy the following properties:
$$\begin{array}{ll} 
& (1)~\{ 2^j\mid j\!\geq\! q_n\}\!\subseteq\!\omega\!\setminus\!\theta_n[\omega ]\cr
& (2)~\theta_n\mbox{ is injective}\cr    
& (3)~\omega\!\setminus\! S_n\! =\!\{ k\!\in\!\omega\mid\theta_n(k)\! =\! k\}\cr    
& (4)~\forall s\!\in\!\omega^{<\omega}\mbox{ strictly decreasing with~}2\!\leq\!\vert s\vert\!\leq\! L~~
\theta_{s^*}[\omega ]\cap S_{s(0)}\mbox{ is infinite}\cr   
& (5)~\forall s\!\in\!\omega^{<\omega}\mbox{ strictly decreasing with~}\vert s\vert\! =\! L\! +\! 1~~
\theta_{s^*}[\omega ]\cap S_{s(0)}\! =\!\emptyset\cr   
& (6)~\theta_n\big( 2^{q_n}\!\cdot\! (j\! +\! 1)\big)\! -\!\theta_n(2^{q_n}\!\cdot\! j)\!\leq\! 
2^{q_n+1}\!\cdot\! 3
\end{array}$$ 
\end{lem}

\noindent {\bf Proof.} (1) If $k\!\in\! S_n$, then $\theta_n(k)$ is a multiple of an odd number $\geq\! 3$. It remains to note that $\{ 2^j\mid j\!\geq\! q_n\}\!\subseteq\! S_n$ and $\theta_n$ is the identity on 
$\omega\!\setminus\! S_n$.\bigskip

\noindent (2) If $L\! =\! 1$, then $\theta_n$ is a bijection between $\omega\!\setminus\! S_n$ onto itself on one side, and from $S_n$ onto 
$$\{ 2^{q_n}\!\cdot\! (2\!\cdot\! k\! +\! 1)\mid k\!\geq\! 1\}\!\subseteq\! S_n$$ 
on the other side. So we may assume that $L\!\geq\! 2$.

\vfill\eject

 The map $\theta'\! :\!\{ j\!\in\!\omega\mid j\!\geq\! 1\}\!\rightarrow\!\{ 3k\mid k\!\geq\! 1\}$ defined by $\theta'(j)\! :=\! 3j$ is a bijection. The difference between $\theta$ and $\theta'$ is that $\theta$ exchanges $\theta'(2^{31}\!\cdot\! 3\!\cdot\! k)$ and 
$\theta'(2^{31}\!\cdot\! 3\!\cdot\! k\! +\! 1)$ for each $k\!\geq\! 1$ with $k\!\notin\! P_L$. Thus $\theta$ is a bijection too. This implies that $\theta_n$ is a bijection between $\omega\!\setminus\! S_n$ onto itself on one side, and from $S_n$ onto 
$\{ 2^{q_n}\!\cdot\! 3\!\cdot\! k\mid k\!\geq\! 1\}\!\subseteq\! S_n$ on the other side.\bigskip

\noindent (3) It is enough to check that $\theta$ is fixed point free. We may assume that $L\!\geq\! 2$. In the first two cases of the definition of $\theta$, $\theta (j)\!\geq\! 3j\! >\! j$. In the last case, if $\theta (j)\! =\! j$, then $2j\! =\! 3$, which is absurd.\bigskip

\noindent (4) Fix $r\!\in\!\omega$. Note that 
$$\theta_{s^*}(2^{q_{s(0)}+r})\! =\!\theta_{<s(1),\cdots ,s(\vert s\vert\! -\! 2)>}(2^{q_{s(0)}+r}\!\cdot\! 3)\! =\!\cdots\! =\!\theta_{s(1)}(2^{q_{s(0)}+r}\!\cdot\! 3^{\vert s\vert -2})\! =\! 
2^{q_{s(0)}+r}\!\cdot\! 3^{\vert s\vert -1}\!\in\! S_{s(0)}.$$
(5) We may assume that $L\!\geq\! 2$. Note first that if $\theta (3j)\! =\! 2^{31}\!\cdot\! 3\!\cdot\! m$, then 
$j\!\in\! P_L$. Indeed, we argue by contradiction. If $3j\!\notin\!\big( M_L\cup (M_L\! +\! 1)\big)$, then 
$\theta (3j)\! =\! 3^2\!\cdot\! j$, so that $3j\! =\! 2^{31}\!\cdot\! m$. In particular, $m$ is a multiple of $3$, which implies that $m$ is of the form $2^p\!\cdot\! 3^{l+1}$ with $l\! <\! L\! -\! 2$ and $j\! =\! 2^{31+p}\!\cdot\! 3^l$. If 
$3j\! =\! 2^{31}\!\cdot\! 3\!\cdot\! k$ and $k\!\notin\! P_L$, then $\theta (3j)\! =\! 3^2\!\cdot\! j\! +\! 3$, so that 
$3j$ is a multiple of $2^{31}$ and $3j\! +\! 1\! =\! 2^{31}\!\cdot\! m$, which is absurd. Finally, 
$3j$ cannot be of the form $2^{31}\!\cdot\! 3\!\cdot\! k\! +\! 1$.\bigskip

 If $\theta_{s^*}(J)\!\in\! S_{s(0)}$, then $\theta_{s(\vert s\vert -1)}(J)\!\in\! S_{s(\vert s\vert -2)}$ is a multiple of $3$ and is of the form $2^{q_{s(\vert s\vert -2)}}\!\cdot\! 3\!\cdot\! j$. Thus 
 $\theta_{<s(\vert s\vert -2),s(\vert s\vert -1)>}(J)\! =\!
\theta_{s(\vert s\vert -2)}(2^{q_{s(\vert s\vert -2)}}\!\cdot\! 3\!\cdot\! j)\! =\! 
2^{q_{s(\vert s\vert -2)}}\!\cdot\!\theta (3j)\!\in\! S_{s(\vert s\vert -3)}$ is of the form 
$2^{q_{s(\vert s\vert -3)}}\!\cdot\! k$. Thus $\theta (3j)\! =\! 2^{q_{s(\vert s\vert -3)}-q_{s(\vert s\vert -2)}}\!\cdot\! k$ is of the form $2^{31}\!\cdot\! 3\!\cdot\! m$ since 
$$q_{s(\vert s\vert -3)}-q_{s(\vert s\vert -2)}\!\geq\! q_{s(\vert s\vert -2)+1}-q_{s(\vert s\vert -2)}\!\geq\! 32^{q_1}\! -\! q_1\! =\! q_2-q_1\! =\! 31.$$ 
The previous point implies that $j\!\in\! P_L$. Thus 
$$\theta_{s^*}(J)\! =\!\theta_{<s(1),\cdots ,s(\vert s\vert\! -\! 2)>}(2^{q_{s(\vert s\vert -2)}+p}\!\cdot\! 3^{l+1})
\! =\!\theta_{<s(1),\cdots ,s(N)>}(2^{q_{s(\vert s\vert -2)}+p}\!\cdot\! 3^{\vert s\vert +l-N-1})$$
as in (4), as long as $\vert s\vert\! +\! l\! -\! N\! -\! 1\!\leq\! L\! -\! 2$, i.e., $l\! +\! 2\!\leq\! N$. In other words, 
$$\theta_{s^*}(J)\! =\!\theta_{<s(1),\cdots ,s(l+2)>}(2^{q_{s(\vert s\vert -2)}+p}\!\cdot\! 3^{L-2})\! =\!
\theta_{<s(1),\cdots ,s(l+1)>}(2^{q_{s(\vert s\vert -2)}+p}\!\cdot\! 3^{L-1}).$$ 
As $\theta_{s^*}(J)\!\in\! S_{s(0)}$, $q_{s(\vert s\vert -2)}+p\!\geq\! q_{s(0)}$. Now  
$\theta_{s(l+1)}(2^{q_{s(\vert s\vert -2)}+p}\!\cdot\! 3^{L-1})\! =\! 
2^{q_{s(\vert s\vert -2)}+p}\!\cdot\! 3^L\! +\! 3$ is not in $S_{s(l)}$, which is the desired contradiction.\bigskip

\noindent (6) We may assume that $L\!\geq\! 2$. Note that $\theta (j\! +\! 1)\! -\!\theta (j)\!\in\!\{ 3,6,-3\}$, so that 
$$\theta_n\big( 2^{q_n}\!\cdot\! (j\! +\! 1)\big)\! -\!\theta_n(2^{q_n}\!\cdot\! j)\!\leq\! 
2^{q_n+1}\!\cdot\! 3.$$
This finishes the proof.\hfill{$\square$}\bigskip

\noindent {\bf Notation.}  We define, for $s\!\in\!\omega^{<\omega}\!\setminus\!\{\emptyset\}$, 
$s^-\! :=<s(0),\cdots ,s(\vert s\vert -2)>$ and 
$$s^{-1}\! :=<s(\vert s\vert -1),\cdots ,s(0)>.$$
We then set, for $n\!\in\!\omega$,  $\mathbb{D}_n^1\! :=\! N_{t_n0}$. We extend the definition of 
$g_n\! :\! N_{t_n0}\!\rightarrow\! N_{t_n1}$ after Theorem \ref{infi} to any ${L\!\geq\! 2}$, writing $g^L_n$ instead of $g_n$ when a confusion is possible. Fix $L\!\geq\! 2$. We set, for $n\!\in\!\omega$,  
$$\mathbb{D}^L_n\! :=\!\big\{\alpha\!\in\! N_{t_n0}\mid\forall m\!\leq\! n~~\alpha (2^{q_n}\!\cdot\! 3^m)\!\not=\!
\alpha (2^{q_n}\!\cdot\! 3^{m+1})\big\}$$ 
and $\mathbb{G}_L\! :=\!\bigcup_{n\in\omega}~\mbox{Graph}({g^L_n}_{\vert\mathbb{D}^L_n})$, so that 
$\mathbb{D}^L_n$ is a nonempty clopen subset of $2^\omega$ and $\mathbb{G}_L$ is a $\boratwo$ digraph on 
$2^\omega$.

\begin{lem} \label{propf++} Fix $L\!\geq\! 1$.\smallskip

(a) $\big( 2^\omega ,({g^L_n}_{\vert\mathbb{D}^L_n})\big)$ is a strongly complex situation.\smallskip

(b) $g_s(\alpha )\!\not=\! g_{s^-}(\alpha )$ if $s\!\in\!\omega^{<\omega}$ is strictly increasing, 
$2\!\leq\!\vert s\vert\!\leq\! L$, $\alpha\!\in\!\mathbb{D}^L_{s(\vert s\vert -1)}$ and $g_s(\alpha ),g_{s^-}(\alpha )$ are defined.\smallskip

(c) $g_s(\alpha )\! =\! g_{s^-}(\alpha )$ if $s\!\in\!\omega^{<\omega}$ is strictly increasing, 
$\vert s\vert\! =\! L\! +\! 1$ and $g_s(\alpha ),g_{s^-}(\alpha )$ are defined.\end{lem}

\noindent {\bf Proof.} (a) It is known that $2^\omega$ is a nonempty zero-dimensional perfect Polish space. Note that $g_n$ is defined, partial, and continuous since $g_n(\alpha )(k)$ depends only on $\alpha\vert\big(\theta_n(k)\! +\! 1\big)$. If $x\!\in\! 2^{<\omega}$, then by Lemma \ref{deff++}.(2) 
$g_n[N_{t_n0x}]\! =\!\big\{\beta\!\in\! N_{t_n1}\mid\forall 2^{q_n}\! <\! i\! <\! 2^{q_n}\! +\! 1\! +\!\vert x\vert ~~
\big( i\! =\!\theta_n(k)\Rightarrow\beta (k)\! =\! (t_n0x)(i)\big)\big\}$ is clopen, so that $g_n$ is open and onto with clopen domain and range, and 
${g^L_n}_{\vert\mathbb{D}^L_n}$ is partial continuous open with clopen domain and range. If $\alpha\!\in\! 2^\omega$ and $l\!\in\!\omega$, then there is $n$ with $\alpha\vert l\! =\!\psi (n)$. Pick $\alpha_0\!\in\!\mathbb{D}^L_n$. Then $\big(\alpha_0,g^L_n(\alpha_0)\big)\!\in\!\mathbb{G}_L\cap N_{\alpha\vert l}^2$. Finally, assume that $U$ is an open subset of $2^\omega$ and $\alpha\!\in\! U\cap\mathbb{D}^L_n$. Then we can find $x\!\in\! 2^{<\omega}$ such that 
$\alpha\!\in\! N_{t_n0x}\!\subseteq\! U\cap\mathbb{D}^L_n$. If $\beta (k)\! =\alpha (k)$ for $k\!\in\!\theta_n[\omega ]\cup 2^{q_n+1+\vert x\vert}$, then 
$\beta\!\in\! N_{t_n0x}\!\subseteq\! U\cap\mathbb{D}^L_n$ and $g_n(\beta )\! =\! g_n(\alpha )$. As $\omega\!\setminus\!\theta_n[\omega ]$ is infinite by Lemma 
\ref{deff++}.(1), $\omega\!\setminus\! (\theta_n[\omega ]\cup 2^{q_n+1+\vert x\vert})$ is also infinite, and the set of such $\beta$'s is not countable.\bigskip

\noindent (b) We argue by contradiction. We set $k\! :=\! 2^{q_{s(\vert s\vert -1)}}$. As $(q_n)$ and $s$ are strictly increasing and $\vert s\vert\!\geq\! 2$, 
$k\!\not=\! 2^{q_{s(0)}}$, and $g_{s(0)}(\beta )(k)\! =\!\beta\big(\theta_{s(0)}(k)\big)$ if $\beta\!\in\! N_{t_{s(0)}0}$. Thus 
$$g_{s(0)}(\beta )(k)\! =\!\beta\big(\theta_{s(0)}(2^{q_{s(0)}}\!\cdot\! 2^{q_{s(\vert s\vert -1)}-q_{s(0)}})\big)\! =\!
\beta\big( 2^{q_{s(0)}}\!\cdot\!\theta (2^{q_{s(\vert s\vert -1)}-q_{s(0)}})\big)\! =\!
\beta (2^{q_{s(\vert s\vert -1)}}\!\cdot\! 3)$$
if $\beta\!\in\! N_{t_{s(0)}0}$. As $2^{q_{s(\vert s\vert -1)}}\!\cdot\! 3\!\not=\! 2^{q_{s(1)}}$,  
$$\begin{array}{ll}
g_{<s(0),s(1)>}(\beta )(k)\!\!\!\!\!
& =\! g_{s(1)}(\beta )(2^{q_{s(\vert s\vert -1)}}\!\cdot\! 3)\! =\!
\beta\big(\theta_{s(1)}(2^{q_{s(\vert s\vert -1)}}\!\cdot\! 3)\big)\cr 
& =\!\beta\big( 2^{q_{s(1)}}\!\cdot\!\theta (2^{q_{s(\vert s\vert -1)}-q_{s(1)}}\!\cdot\! 3)\big)\! =\!
\beta (2^{q_{s(\vert s\vert -1)}}\!\cdot\! 3^2)
\end{array}$$ 
if $\beta\!\in\! N_{t_{s(1)}0}$ and $g_{s(1)}(\beta )\!\in\! N_{t_{s(0)}0}$. Similarly, 
$g_s(\alpha )(k)\! =\!\alpha (2^{q_{s(\vert s\vert -1)}}\!\cdot\! 3^{\vert s\vert})$ since $2\!\leq\!\vert s\vert\!\leq\! L$, which contradicts the fact that $\alpha\!\in\!\mathbb{D}^L_{s(\vert s\vert -1)}$ since 
$\vert s\vert\! -\! 1\!\leq\! s(\vert s\vert\! -\! 1)$ cause $s$ is stricly increasing.\bigskip

\noindent (c) As above, $g_s(\alpha )(k)\! =\!\alpha\big(\theta_{s^{-1}}(k)\big)$ if $k\!\not=\! 2^{q_{s(0)}}$. As 
$\theta_{(s^-)^{-1}}(k)$ is in $\theta_{(s^-)^{-1}}[\omega ]$, it is not in $S_{s(\vert s\vert -1)}$ by Lemma \ref{deff++}.(5). Thus 
$\theta_{s^{-1}}(k)\! =\!\theta_{(s^-)^{-1}}(k)$ by Lemma \ref{deff++}.(3) and we are done.\hfill{$\square$}\bigskip

 We introduce some finitary objects used in the construction of our homomorphism. Fix $L\!\geq\! 1$.\bigskip

\noindent\bf Notation.\rm\ We inductively define, for each $l\!\in\!\omega$,\smallskip

\noindent - a subset $X_l$ of $2^{\leq l}$ such that $2^\omega$ is the disjoint union of $(N_x)_{x\in X_l}$,\smallskip

\noindent - an oriented graph $B_l$ on $X_l$, containing finite approximations of $\mathbb{G}_L$ and providing some control on the cycles,\smallskip

\noindent - a map $\varphi_l\! :\! B_{l+1}\!\rightarrow\! l$ such that $t_{\varphi_l(y,x)}0\!\subseteq\! y$ and $t_{\varphi_l(y,x)}1\!\subseteq\! x$ if $(y,x)\!\in\! B_{l+1}$, giving the number of the function approximated by an element of $B_{l+1}$,\smallskip

\noindent - an uogas $A_l$ on $X_l$ contained in $B_l$, containing finite approximations of a subset of 
$\mathbb{G}_L$ with acyclic symmetrization,\smallskip

\noindent - a subset $E_l$ of $X_l$, determining the future elements of $X_{l+1}$.\smallskip

 We first perform the construction of these objects for $l\! =\! 0$, sometimes more.\smallskip

\noindent - We set $X_0\! :=\!\{\emptyset\}$, so that $X_0\!\subseteq\! 2^{\leq 0}$ and $2^\omega$ is the disjoint union of $(N_x)_{x\in X_0}$.\smallskip

\noindent - Generally speaking, we set 
$$B_l\! :=\!\big\{ (y,x)\!\in\! X_l^2\mid(N_y\!\times\! N_x)\cap\mathbb{G}_L\!\not=\!\emptyset\wedge
\exists i\! <\!\vert y\vert ,\vert x\vert ~~y(i)\!\not=\! x(i)\big\}\mbox{,}$$ 
so that $B_0\! =\!\emptyset$. Note that $B_l$ is an oriented graph on $X_l$ since 
$B_l\!\subseteq <_{\mbox{lex}}$.\bigskip

\noindent - We set $A_0\! :=\!\emptyset$.\smallskip

\noindent - If $x\!\in\! X_l$, then the decision of putting $x$ in $E_l$ or not is made by induction on $M_x$, related to $A_l$, defined before Lemma \ref{disj}. If $M_x\! =\! 1$, then we put $x$ in $E_l$. We set, when the decision of putting $y$ in $E_l$ is made,
$$l_y\! :=\!\left\{\!\!\!\!\!\!\!\!
\begin{array}{ll}
& \vert y\vert\mbox{ if }y\!\notin\! E_l\mbox{,}\cr
& \vert y\vert\! +\! 1\mbox{ if }y\!\in\! E_l.
\end{array}
\right.$$
Assume that the decision of putting $x$ in $E_l$ is made if $1\!\leq\! M_x\!\leq\! m$, which is the case for 
$m\! =\! 1$, and fix $x\!\in\! X_l$ with $M_x\! =\! m\! +\! 1$, so that $x\!\notin\!\mbox{min}_{X_l}$ and 
$\mbox{Pred}(x)\!\not=\!\emptyset$. Fix $y\!\in\!\mbox{Pred}(x)$, so that $(y,x)\!\in\! A_l$. We put 
$n\! :=\!\varphi_{l-1}(y,x)$. We say that $x$ is {\bf $y$-expandable} if $\theta_n(\vert x\vert )\! <\! l_y$. Generally speaking, we put $x$ in $E_l$ exactly when $x$ is $y$-expandable for each $y\!\in\!\mbox{Pred}(x)$ and, when there is $q\!\in\!\omega$ with $t_q\!\in\! E_l$, there is no $1\!\leq\! i\! <\!\vert p_{t_q}\vert$ with 
$x\! =\! p_{t_q}(i)$. In particular, $E_0\! =\!\{\emptyset\}$.\bigskip

 Assume that our objects have been constructed for $p\!\leq\! l$, which is the case for $l\! =\! 0$.\smallskip

\noindent - We set 
$X_{l+1}\! :=\!\big\{ s\varepsilon\mid s\!\in\! E_l\wedge\varepsilon\!\in\! 2\big\}\cup (X_l\!\setminus\! E_l)$, so that $X_{l+1}\!\subseteq\! 2^{\leq l+1}$ and $2^\omega$ is the disjoint union of $(N_x)_{x\in X_{l+1}}$.\smallskip

\noindent - If $(y,x)\!\in\! B_{l+1}$, then there is a unique $n\!\in\!\omega$ with 
$(N_y\!\times\! N_x)\cap\mbox{Graph}(g^L_n)\!\not=\!\emptyset$, and $2^{q_n}\! =\!\vert y\wedge x\vert$. Note that $n\!\leq\! q_n\! <\! 2^{q_n}\! <\!\vert y\vert\!\leq\! l\! +\! 1$. We set $\varphi_l(y,x)\! :=\! n$, which defines $\varphi_l\! :\! B_{l+1}\!\rightarrow\! l$.\smallskip

\noindent - We now define $A_{l+1}$. If $t_q\!\in\!\mbox{max}_{X_l}\cap E_l$, then we put $(t_q0,t_q1)$
in $A_{l+1}$. If $y\!\in\! X_l\!\setminus\!\mbox{max}_{X_l}$, then let $x$ be the unique element of 
$\mbox{Succ}(y)$, so that $(y,x)\!\in\! A_l$. As $A_l\!\subseteq\! B_l$, we can define 
$n\! :=\!\varphi_{l-1}(y,x)$. If $x\!\notin\! E_l$, then we put 
$$\left\{\!\!\!\!\!\!\!
\begin{array}{ll}
& (y,x)\mbox{ if }y\!\notin\! E_l\mbox{,}\cr
& (y\eta ,x)\mbox{ if }y\!\in\! E_l^-\! :=\! E_l\!\setminus\!\{ t_m\mid m\!\in\!\omega\}\wedge\eta\!\in\! 2
\mbox{,}\cr
& (y1,x),(y0,y1)\mbox{ if }y\! =\! t_q\!\in\! E_l
\end{array}
\right.$$
in $A_{l+1}$. If $x\!\in\! E_l$, then $y\!\notin\!\{ t_m\mid m\!\in\!\omega\}$, $x$ is $y$-expandable and 
$\theta_n(\vert x\vert )\! <\! l_y$. We put $(y\eta ,x\varepsilon )$ in $A_{l+1}$ if $\eta\!\in\! 2^{\leq 1}$, 
$y\eta\!\in\! X_{l+1}$ and $\varepsilon\! =\! (y\eta )\big(\theta_n(\vert x\vert )\big)$. We first check the following announced facts.
 
\begin{lem} \label{uogas} Let $l\!\in\!\omega$. Then $A_l$ is an uogas on $X_l$ contained in $B_l$.\end{lem} 

\noindent\bf Proof.\rm\ We argue by induction on $l$, and the case $l\! =\! 0$ is clear. Note that 
$A_{l+1}\!\subseteq\! B_{l+1}$, by definition of $A_{l+1}$, so that $A_{l+1}$ is an oriented graph on 
$X_{l+1}$. By  definition, $A_{l+1}$ is unambiguously oriented. We argue by contradiction to see that $s(A_{l+1})$ is acyclic. Let $(u_j)_{j\leq N}$ be a $s(A_{l+1})$-cycle. We choose ${j\!\leq\! N}$ such that $u_j$ is $<_{\mbox{lex}}$-minimal, and we may assume that $0\! <\! j\! <\! N$. Then 
$(u_j,u_{j-1}),(u_j,u_{j+1})\!\in\! A_{l+1}$, so that $u_{j-1}\! =\! u_{j+1}$ since $A_{l+1}$ is unambiguously oriented, which is absurd.\hfill{$\square$}

\begin{lem} \label{injec} Let $l\!\in\!\omega$, and $x\!\in\! X_{l+1}$. Then 
$\vert p_x\vert\!\leq\! l\! +\! 1$.\end{lem} 
 
\noindent\bf Proof.\rm\ We argue by induction on $l$. For $l\! =\! 0$, $X_1\! =\!\{ (0),(1)\}$, 
$p_{(1)}\! =\!\big( (1)\big)$ and $p_{(0)}\! =\!\big( (0)\big)$, so we are done. Let $x\!\in\! X_{l+2}$. Note first that if $i\! <\! L\! :=\!\vert p_x\vert$, then we can find a unique couple $(x_i,\varepsilon_i)$ in 
$X_{l+1}\!\times\! 2^{\leq 1}$ with $p_x(i)\! =\! x_i\varepsilon_i$. Let $i\! <\! L\! -\! 1$. By definition of 
$A_{l+2}$, either $(x_i,x_{i+1})\!\in\! A_{l+1}$, or there is $q\!\in\!\omega$ such that $t_q\!\in\! E_{l+1}$ and 
$\big( p_x(i),p_x(i\! +\! 1)\big)\! =\! (t_q0,t_q1)$. In particular, $(x_i)_{i<L}$ is $\leq_{\mbox{lex}}$-increasing. By definition of $E_{l+1}$, there is at most one $i\! <\! L\! -\! 1$ for which there is $q\!\in\!\omega$ such that $t_q\!\in\! E_{l+1}$ and $\big( p_x(i),p_x(i\! +\! 1)\big)\! =\! (t_q0,t_q1)$. If such a $q$ does not exist, then 
$(x_i)_{i<L}\! =\! p^{s(A_{l+1})}_{x_0,x_{L-1}}$ and $L\!\leq\! l\! +\! 1$, by induction assumption. If now there is such a $q$, then $(x_0,\cdots ,x_i,x_{i+2},\cdots ,x_{L-1})\! =\! p^{s(A_{l+1})}_{x_0,x_{L-1}}$ and 
$L\! -\! 1\!\leq\! l\! +\! 1$, by induction assumption, so that $L\!\leq\! l\! +\! 2$.\hfill{$\square$}

\begin{lem} \label{inc} Fix $L\!\geq\! 1$. Then there is $(L_n)_{n\in\omega}\!\in\! (\omega\!\setminus\!\{ 0\} )^\omega$ strictly increasing satisfying the following properties:\smallskip

(a) $t_n\!\in\! E_{L_n}\!\setminus\! (\bigcup_{k<L_n}~E_k)$,\smallskip
 
(b) $X_{L_n}\cap X_{L_{n+1}}\! =\!\emptyset$,\smallskip
 
(c) $2^{q_n}\!\leq\! L_n\! <\! L_n\! +\! (6L_n)^{L_n-1}\!\leq\! 2^{q_{n+1}}$.\end{lem} 
 
\noindent\bf Proof.\rm\ We construct $L_n$ inductively on $n$. We set $L_0\! :=\! 1$, which is correct since 
$t_0\! =\! (0)$ is in $E_1\!\setminus\! E_0$, $p_0\! =\! 0$ and $p_1\! =\! 1$. As $X_1\! =\!\{ (0),(1)\}$, 
$\vert x\vert\! >\! 0$ if $x\!\in\! X_{L_0}$. Assume that $L_n$ has been constructed.\bigskip

\noindent\bf Claim.\it\ Let $l\!\geq\! L_n\! +\! (3\!\cdot\! 2^{q_n+1})^{L_n-1}$. Then $X_l$ does not meet $X_{L_n}$. In fact, we can add $m\!\geq\! 1$ coordinates at least to $x\!\in\! X_{L_n}$ in $X_l$ if 
$l\!\geq\! L_n\! +\! m\!\cdot\! (3\!\cdot\! 2^{q_n+1})^{L_n-1}$ and no $t_q$ appears in some $E_k$ with 
$L_n\! <\! k\! <\! l$.\rm\bigskip

 Indeed, fix $x\!\in\! X_{L_n}$. We want to extend properly $x$ in some $X_l$ with $l\! >\! L_n$, and to give an estimate on $l$. We proceed by induction on $M_x$, which was defined just before Lemma \ref{disj}. If $M_x\! =\! 1$, then $l\! :=\! L_n\! +\! 1\!\leq\! L_n\! +\! (3\!\cdot\! 2^{q_n+1})^{L_n-1}$ is suitable.  Assume now that $M_x\! =\! 2$, which gives $y\!\in\!\mbox{min}_{X_{L_n}}$ with $x\!\in\! p_y$, and in fact $x\! =\! p_y(1)$. In particular, 
$(y,x)\!\in\! A_{L_n}\!\subseteq\! B_{L_n}$ and $N\! :=\!\varphi_{L_n-1}(y,x)\! <\! L_n\! -\! 1$ is defined. The definition of $g^L_N$ and Lemma \ref{deff++}.(6) show that if $\alpha\!\in\! N_y$, then we can find 
$l\!\leq\! L_n\! +\! 3\!\cdot\! 2^{q_N+1}$ and $z\!\subseteq\!\alpha$ with $z\!\in\! X_l$ and 
${z(\vert z\vert\! -\! 1)\! =\! g^L_N(\alpha )(\vert x\vert )}$. Note that $q_N\!\leq\! q_n$ since 
$(q_n)$ is increasing, $L_n\!\geq\! 2$ since $B_0\! =\! B_1\! =\!\emptyset$, and 
$L_n\! <\! L_n\! +\! 3\!\cdot\! 2^{q_n+1}\!\leq\! L_n\! +\! 6L_n\!\leq\! L_n\! +\! (6L_n)^{L_n-1}\!\leq\! 
2^{q_{n+1}}$. Thus $x$ is properly extended in $X_{L_n+3\cdot 2^{q_n+1}}$, and also in 
$X_{L_n+(3\cdot 2^{q_n+1})^{L_n-1}}$. This argument also shows that $x$ can be extended by $m$ coordinates in $X_{L_n+ m\cdot 3\cdot 2^{q_n+1}}$ if 
$L_n\! <\! L_n\! +\! m\!\cdot\! 3\!\cdot\! 2^{q_n+1}\!\leq\! 2^{q_{n+1}}$. More generally, this argument shows that $x$ can be properly extended in $X_{L_n+(3\cdot 2^{q_n+1})^{L_n-1}}$, by Lemma \ref{injec}, which implies that 
$M_x\!\leq\! L_n$.\hfill{$\diamond$}\bigskip

 Now note that $\vert\psi (n\! +\! 1)\vert\!\leq\! n\! +\! 1\!\leq\! 2^n\!\leq\! 2^{q_n}\!\leq\! L_n\! <\! 
L_n\! +\! (3\!\cdot\! 2^{q_n+1})^{L_n-1}\!\leq\! 2^{q_{n+1}}$. We will extend $\psi (n\! +\! 1)$ until we reach 
$t_{n+1}\! =\!\psi (n\! +\! 1)0^{2^{q_{n+1}}-\vert\psi (n+1)\vert}$ in some $X_q$ with $q\!\geq\! 2^{q_{n+1}}$. Note first that there is $x\!\in\! X_{L_n}$ with $x\!\subseteq\!\psi (n\! +\! 1)0^\infty$, and $\psi (n\! +\! 1)\!\subseteq\! x$ since 
$\vert x\vert\! >\! n$.\bigskip

 In order to do this, we add $m\! :=\! 2^{p_{n+1}}\! -\!\vert x\vert$ coordinates to $x$. By the claim, this will be possible in $X_q$ for some $2^{q_{n+1}}\!\leq\! q\!\leq\! L_n\! +\! m\!\cdot\! (3\!\cdot\! 2^{q_n+1})^{L_n-1}$ if no $t_r$ appears in some $E_k$ with $L_n\! <\! k\! <\! q$. This is the case since $q\! <\! L_n\! +\! (m\! +\! 1)\!\cdot\! 
(3\!\cdot\! 2^{q_n+1})^{L_n-1}\!\leq\! 2^{q_{n+1}}\! +\! 2^{q_{n+1}}\!\cdot\! (6L_n)^{L_n-1}\!\leq\! 
(2^{q_{n+1}})^2\! <\! 2^{q_{n+2}}$ because $q_{n+2}\! =\! 32^{q_{n+1}}\! >\! 2q_{n+1}$. It remains to check that 
$L_{n+1}\! +\! (6L_{n+1})^{L_{n+1}-1}\!\leq\! 2^{q_{n+2}}$. As $L_{n+1}\!\leq\! (2^{q_{n+1}})^2$, 
$L_{n+1}\! +\! (6L_{n+1})^{L_{n+1}-1}\!\leq\! (6L_{n+1})^{L_{n+1}}\!\leq\! 
2^{(2q_{n+1}+3)\cdot 2^{2q_{n+1}}}$. It remains to note that  
$(2q_{n+1}\! +\! 3)\!\cdot\! 2^{2q_{n+1}}\!\leq\! 2^{5q_{n+1}}\! =\! 32^{q_{n+1}}\! =\! q_{n+2}$.\hfill{$\square$}

\begin{cor} \label{ext} For each $x\!\in\! 2^{<\omega}$ there is $l\!\in\!\omega$ such that $x\!\in\! X_l$.\end{cor} 
 
\noindent\bf Proof.\rm\ We argue by induction on $\vert x\vert$. For $x\! =\!\emptyset$, we can take 
$l\! :=\! 0$. Assume that $x\!\in\! X_l$ and $\varepsilon\!\in\! 2$. If there is $l'\! >\! l$ with 
$x\varepsilon\!\in\! X_{l'}$, then we are done. Otherwise, $x\!\in\! X_{l'}$ for each $l'\!\geq\! l$. Using Lemma \ref{inc}, we choose $n\!\in\!\omega$ such that $L_n\!\geq\! l$. Then $x\!\in\! X_{L_n}\cap X_{L_{n+1}}$, which is absurd.\hfill{$\square$}

\begin{lem} \label{Bl} We work with $L\! =\! 1$. Let $l\!\in\!\omega$, $(y,x)\!\in\! B_{l+1}$ and 
$p\! :=\! p^{s(A_{l+1})}_{y,x_{C(y)}}$. Then\smallskip

(a) $x\! =\! p(j)$ for some $1\!\leq\! j\! <\!\vert p\vert$,\smallskip

(b) $\varphi_l\big( p(j\! -\! 1),p(j)\big)\! =\!\varphi_l(y,x)\! =\!\mbox{min}_{i<j}~\varphi_l\big( p(i),p(i\! +\! 1)\big)$,\smallskip

(c) $\Big(\varphi_l\big( p(i),p(i\! +\! 1)\big)\Big)_{i<j}$ is injective.\end{lem} 
 
\noindent\bf Proof.\rm\ We argue by induction on $l$. We are done if $l\! =\! 0$ since $B_1\! =\!\emptyset$. So assume that $(y,x)\!\in\! B_{l+2}$, which gives $n\! :=\!\varphi_{l+1}(y,x)\!\leq\! l$ with 
$t_n0\!\subseteq\! y$ and $t_n1\!\subseteq\! x$. In particular, $l\! +\! 1\!\geq\! L_n$. If $l\! +\! 1\! =\! L_n$, then 
$(y,x)\! =\! (t_n0,t_n1)$ and $p\! =\! (t_n0)p^{s(A_{l+2})}_{t_n1,x_{C(t_n0)}}$ by Lemma \ref{am1} and since $A_{l+2}\!\subseteq <_{\mbox{lex}}$. Thus $j\! =\! 1$ is convenient. So we may assume that $l\!\geq\! L_n$. Note that we can find $(y',x')\!\in\! B_{l+1}$ with $t_n0\!\subseteq\! y'\!\subseteq\! y$ and 
$t_n1\!\subseteq\! x'\!\subseteq\! x$. By the induction assumption, $p'\! :=\! p^{s(A_{l+1})}_{y',x_{C(y')}}$ is defined and there is $1\!\leq\! j'\! <\!\vert p'\vert$ with $x'\! =\! p'(j')$.\bigskip

\noindent\bf Case 1.\rm\ We cannot find $q\!\in\!\omega$ with $t_q\!\in\! p^{s(A_{l+1})}_{y',x'}\cap E_{l+1}$.\bigskip

 Let $\varepsilon_0\!\in\! 2^{\leq 1}$ such that $y\! =\! y'\varepsilon_0$. Note that there is a unique 
$\varepsilon_1\!\in\! 2^{\leq 1}$ such that $\big( y,p'(1)\varepsilon_1\big)$ is in $A_{l+2}$, by definition of $A_{l+2}$. Similarly, if $1\!\leq\! i\! <\! j'$, then there is a unique $\varepsilon_{i+1}\!\in\! 2^{\leq 1}$ such that $\big( p'(i)\varepsilon_i,p'(i\! +\! 1)\varepsilon_{i+1}\big)\!\in\! A_{l+2}$. Note that 
$p^{s(A_{l+2})}_{y,p'(j')\varepsilon_{j'}}$ is defined and equal to $\big( p'(i)\varepsilon_i\big)_{i\leq j'}$. We will be able to set $j\! :=\! j'$ if we prove that 
$x\! =\! x'\varepsilon_{j'}$,  by Lemma \ref{am1} and since $A_{l+2}\!\subseteq <_{\mbox{lex}}$ again. This is the case if $\varepsilon_{j'}\! =\!\emptyset$, so we may assume that $\varepsilon_{j'}\!\in\! 2$, which implies that $x'\!\in\! E_{l+1}$. Let $\varepsilon\! :=\! x(\vert x'\vert )$, so that $x\! =\! x'\varepsilon$ and we have to see that $\varepsilon_{j'}\! =\!\varepsilon$. Note that 
$p'(i\! +\! 1)\varepsilon_{i+1}$ is $\big( p'(i)\varepsilon_i\big)$-expandable if $\varepsilon_{i+1}\!\in\! 2$, for each $i\! <\! j'$. More precisely, if 
$n_i\! :=\!\varphi_l\big( p'(i),p'(i\! +\! 1)\big)$, then 
$\varepsilon_{i+1}\! =\!\big( p'(i)\varepsilon_i\big)\Big(\theta_{n_i}\big(\vert p'(i\! +\! 1)\vert\big)\Big)$. This implies that $\varepsilon_{j'}\! =\! y\big(\theta_n(\vert x'\vert )\big)$ since $n_{j'-1}\! =\! n\! =\!\mbox{min}_{i<j'}~n_i$ and 
$(n_i)_{i<j'}$ is injective, by induction assumption. In particular, $y$ is long enough to ensure that 
$\varepsilon_{j'}\! =\!\varepsilon$ since $(y,x)\!\in\! B_{l+2}$ with witness $n$.\bigskip

\noindent\bf Case 2.\rm\ $t_q\!\in\! p^{s(A_{l+1})}_{y',x'}\cap E_{l+1}$ for some $q\!\in\!\omega$, which implies that $l\! +\! 1\! =\! L_q$.\bigskip

 Fix $i_0\!\leq\! j'$ with $p'(i_0)\! =\! t_q$. The definition of $\varepsilon_i$ is as in Case 1 if $i\!\leq\! i_0$. If $\varepsilon_{i_0}\! =\! 0$, then we set 
${\varepsilon_i\! :=\!\emptyset}$ if $i_0\! <\! i\!\leq\! j'$ and we note that 
$p^{s(A_{l+2})}_{y,x}\! =\!\big( p'(0)\varepsilon_0,\cdots ,p'(i_0)0,p'(i_0)1,p'(i_0\! +\! 1),\cdots, x'\big)$ if ${i_0\! <\! j'}$, and 
$p^{s(A_{l+2})}_{y,x}\! =\!\big( p'(0)\varepsilon_0,\cdots ,p'(i_0)0\big)$ if $i_0\! =\! j'$. Apart from that, we argue as in Case 1. If ${\varepsilon_{i_0}\! =\! 1}$, then we set $\varepsilon_i\! :=\!\emptyset$ if $i_0\! <\! i\!\leq\! j'$ and we note that 
$$p^{s(A_{l+2})}_{y,x}\! =\!\big( p'(0)\varepsilon_0,\cdots ,p'(i_0)1,p'(i_0\! +\! 1),\cdots, x'\big)$$ 
if $i_0\! <\! j'$, and $p^{s(A_{l+2})}_{y,x}\! =\!\big( p'(0)\varepsilon_0,\cdots ,p'(i_0)1\big)$ if $i_0\! =\! j'$. Apart from that, we argue as in Case 1.\bigskip

 Note that (b), (c) follow from the previous discussion since 
$\Big(\varphi_l\big( p'(i),p'(i\! +\! 1)\big)\Big)_{i<j'}$ is equal to 
$\Big(\varphi_{l+1}\big( p(i),p(i\! +\! 1)\big)\Big)_{i<j}$, except when $i_0\! <\! j'$ in the first subcase of the Case 2 where one number bigger than the others has been added strictly before the last position.\hfill{$\square$}

\begin{lem} \label{lkyk} Fix $L\!\geq\! 1$. If $n\!\in\!\omega$, $y_{\vert t_n\vert +1}\! :=\! t_n1$ and 
$k\! >\!\vert t_n\vert\! +\! 1$, then for each $\alpha\!\in\! N_{t_n0}$ we can find $l_k\!\in\!\omega$ and 
$y_k\!\in\! X_{l_k}$ such that $(\alpha\vert k,y_k)\!\in\! B_{l_k}$ and $y_{k-1}\!\subseteq\! y_k$; moreover, 
$\{ g_n(\alpha )\}\! =\!\bigcap_{k>\vert t_n\vert}~N_{y_k}$.\end{lem} 
 
\noindent\bf Proof.\rm\ By Corollary \ref{ext}, there is $l_k\!\in\!\omega$ such that 
$\alpha\vert k\!\in\! X_{l_k}$. Note that $l_k\! >\! L_n\! +\! 1$ since $k\! >\!\vert t_n\vert\! +\! 1$. There is 
$y_k\!\in\! X_{l_k}$ such that $g_n(\alpha )\!\in\! N_{y_k}$ since $2^\omega$ is the disjoint union of 
$(N_x)_{x\in X_{l_k}}$. Note that $y_k\!\supseteq\! t_n1$ since $t_n1\!\in\! X_{L_n+1}$ and 
$l_k\! >\! L_n\! +\! 1$, so that $(\alpha\vert k,y_k)\!\in\! B_{l_k}$.\hfill{$\square$}

\section{$\!\!\!\!\!\!$ The main construction}\indent

 We now come to the construction of our homomorphism.
 
\begin{thm} \label{newmain} Let $\big( Z,(f_n)\big)$ be a strongly complex situation satisfying Condition (d) in Theorem \ref{cormain}. Then $(2^\omega,\mathbb{G}_1)\preceq^{\mbox{inj}}_c(Z,A^f)$.\end{thm}

\noindent\bf Proof.\rm\  We construct, inductively on $l$,\bigskip

\noindent $\diamond$ a sequence $(U^l_x)_{x\in X_l}$ of nonempty clopen subsets of $Z$,\smallskip

\noindent $\diamond$ a natural number $\phi (n)$ if $L_n\! <\! l$.\bigskip
 
 We want these objects to satisfy the following conditions, using Lemma \ref{Bl}.
$$\begin{array}{ll} 
& (1)~U^{l+1}_{x\varepsilon}\subseteq U^l_x\mbox{ if }\varepsilon\!\in\! 2^{\leq 1}\wedge 
x\varepsilon\!\in\! X_{l+1}\cr    
& (2)~\mbox{diam}(U^l_x)\!\leq\! 2^{-l}\cr    
& (3)~U^{l+1}_{x0}\cap U^{l+1}_{x1}\! =\!\emptyset\mbox{ if }x0,x1\!\in\! X_{l+1}\cr   
& (4)~U^{l+1}_y\!\subseteq\! D_{\phi (\varphi_l(y,x))}\wedge 
U^{l+1}_x\!\subseteq\! f_{\phi (\varphi_l(y,x))}[U^{l+1}_y]~\mbox{ if }~(y,x)\!\in\! A_{l+1}\cr
& (5)~U^{l+1}_{p(m)}\!\subseteq\! D_{\phi(\varphi_l(y,x))}\mbox{ if }(y,x)\!\in\! B_{l+1}\wedge m\! <\! j~\wedge\cr
& \hfill{\varphi_l\big( p(j\! -\! 1),p(j)\big)\! =\!\varphi_l(y,x)\! =\!
\mbox{min}_{i<j}~\varphi_l\big( p(i),p(i\! +\! 1)\big) }\cr
& (6)~\phi (r)\! >\!\mbox{sup}_{n<r}~\phi (n)
\end{array}$$
Assume that this is done. Fix $\alpha\!\in\! 2^\omega$ and $l\!\in\!\omega$. As $2^\omega$ is the disjoint union of the $N_x$'s for $x\!\in\! X_l$, there is a unique $k_l\!\in\!\omega$ for which $\alpha\vert k_l\!\in\! X_l$. The sequence $(U^l_{\alpha\vert k_l})_{l\in\omega}$ defines $h(\alpha )\!\in\! Z$, using (1) and (2). Note that 
$h\! :\! 2^\omega\!\rightarrow\! Z$ is continuous, and injective by (3) and Corollary \ref{ext} (which implies that 
$k_l$ tends to infinity as $l$ tends to infinity).\bigskip

 Let us prove that $U^{l+1}_y\!\subseteq\! D_{\phi (\varphi_l(y,x))}$ and 
$U^{l+1}_x\!\subseteq\! f_{\phi (\varphi_l(y,x))}[U^{l+1}_y]$ if $(y,x)\!\in\! B_{l+1}$. By Lemma \ref{Bl}, 
$\varphi_l\big( p(j\! -\! 1),p(j)\big)\! =\!\varphi_l(y,x)\! =\!
\mbox{min}_{i<j}~\varphi_l\big( p(i),p(i\! +\! 1)\big)$ and $\varphi_l(y,x)\! <\!\varphi_l\big( p(i),p(i\! +\! 1)\big)$ if 
$i\! <\! j\! -\! 1$. By Lemmas \ref{uogas}, \ref{am1}.(b) and (4), 
$U^{l+1}_{p(i)}\!\subseteq\! D_{\phi (\varphi_l(p(i),p(i+1)))}$ and 
$U^{l+1}_{p(i+1)}\!\subseteq\! f_{\phi (\varphi_l(p(i),p(i+1)))}[U^{l+1}_{p(i)}]$ if $i\! <\! j$. (5) implies that 
$U^{l+1}_{p(m)}\!\subseteq\! D_{\phi (\varphi_l(y,x))}$ if $m\! <\! j$.

\vfill\eject

 We apply Lemma \ref{trans} to 
$V_0\! :=\! U^{l+1}_{p(i)}$, $V_1\! :=\! U^{l+1}_{p(i+1)}$, $m\! :=\!\phi\big(\varphi_l(y,x)\big)$ and 
$$n\! :=\!\phi\Big(\varphi_l\big( p(i),p(i\! +\! 1)\big)\Big)$$ 
when $i\! <\! j\! -\! 1$. We get 
$f_{\phi (\varphi_l(y,x))}[U^{l+1}_{p(i+1)}]\!\subseteq\! f_{\phi (\varphi_l(y,x))}[U^{l+1}_{p(i)}]$, so that 
$$U^{l+1}_x\!\subseteq\! f_{\phi (\varphi_l(y,x))}[U^{l+1}_{p(j-1)}]\!\subseteq\!\cdots\!\subseteq\! 
f_{\phi (\varphi_l(y,x))}[U^{l+1}_y].$$ 
Thus $U^{l+1}_x\!\subseteq\! f_{\phi (\varphi_l(y,x))}[U^{l+1}_y]$, as desired.\bigskip

 Let $(\alpha ,\beta )\!\in\!\mathbb{G}_1$, and $n$ with $\alpha\!\in\! N_{t_n0}$ and $\beta\! =\! g_n(\alpha )$. By Lemma \ref{lkyk}, ${(\alpha\vert k,\beta\vert\vert y_k\vert )\!\in\! B_{l_k}}$ if $k\! >\!\vert t_n\vert\! +\! 1$. By the previous point, $U^{l_k}_{\alpha\vert k}\!\subseteq\! D_{\phi (n)}$ and 
$U^{l_k}_{\beta\vert\vert y_k\vert}\!\subseteq\! f_{\phi (n)}[U^{l_k}_{\alpha\vert k}]$ if 
$k\! >\!\vert t_n\vert\! +\! 1$. In particular, $h(\alpha )\!\in\! D_{\phi (n)}$. As 
$f_{\phi (n)}$ is continuous at $h(\alpha )$, 
$\mbox{diam}(f_{\phi (n)}[U^{l_k}_{\alpha\vert k}])$ converges to zero as $k$ converges to infinity. Thus 
$d\big( f_{\phi (n)}\big( h(\alpha )\big) ,h(\beta )\big)$ is zero and 
$f_{\phi (n)}\big( h(\alpha )\big)\! =\! h(\beta )$.\bigskip

 So it is enough to prove that the construction is possible. We fix a compatible metric with 
$\mbox{diam}(Z)\!\leq\! 1$. We first set $U_\emptyset\! :=\! Z$. Assume that 
$\big( (U_x)_{x\in X_p}\big)_{p\leq l}$ and $\big(\phi (n)\big)_{L_n<l}$ satisfying (1)-(6) have been constructed, which is the case for $l\! =\! 0$.\bigskip

 By Lemmas \ref{uogas} and \ref{am1}.(c), we can set, for each $y\!\in\! X_l$, 
$p^l_y\! :=\! p^{s(A_l)}_{y,x_{C(y)}}$. We choose $u_l\!\in\!\omega^{X_l}$ such that 
$u_l(y)\! :=\!\phi\Big(\varphi_{l-1}\big( y,p^l_y(1)\big)\Big)$ if $y\!\in\! X_l\!\setminus\!\mbox{max}_{X_l}$ (this can be done, by the induction assumption). We define, for $y\!\in\! X_l$, and inductively on $\vert p^l_y\vert$,
$$W_y\! :=\!\left\{\!\!\!\!\!\!\!
\begin{array}{ll} 
& U^l_y\mbox{ if }y\!\in\!\mbox{max}_{X_l}\mbox{,}\cr
& U^l_y\cap f_{u_l(y)}^{-1}(W_x)\mbox{ if }x\!\in\!\mbox{Succ}(y)\mbox{,}
\end{array}
\right.$$
which defines nonempty clopen subsets of $Z$, by (4) of the induction assumption. Note that 
$$f_{u_l(y)}[W_y]\! =\! W_x$$ 
if $x\!\in\!\mbox{Succ}(y)$ since $x\! =\! p^l_y(1)$ and 
$W_x\!\subseteq\! U^l_x\!\subseteq\! f_{\phi (\varphi_{l-1}(y,x))}[U^l_y]$ by the induction assumption.\bigskip

 If $l$ is of the form $L_r$, then Lemma \ref{creation+} applied to $V\! :=\! W_{t_r}$ and 
$m\! :=\!\mbox{sup}_{L_n<l}~\phi (n)\! =\!\mbox{sup}_{n<r}~\phi (n)$ gives $\phi (r)\! >\!\mbox{sup}_{n<r}~\phi (n)$ and nonempty clopen subsets $O_{t_r0},O_{t_r1}$ of $Z$ such that $O_{t_r0}\!\subseteq\! W_{t_r}\cap D_{\phi (r)}$ and 
$O_{t_r1}\!\subseteq\! W_{t_r}\cap f_{\phi (r)}[O_{t_r0}]$.\bigskip

 We will apply, thanks to Lemma \ref{uogas}, Lemma \ref{disj} to 
${\cal T}\! :=\!\big( X_{l+1},A_{l+1},Z,(f_n)\big)$, 
$$d\! :=\!\mbox{max}_{x\in X_{l+1}}~\vert x\vert\mbox{,}$$ 
$u\!\in\!\omega^{X_{l+1}}$ such that $u(y)\! :=\!\phi\Big(\varphi_l\big( y,p^{l+1}_y(1)\big)\Big)$ if 
$y\!\in\! X_{l+1}\!\setminus\!\mbox{max}_{X_{l+1}}$, and $(V_x)_{x\in X_{l+1}}$ defined as follows. If 
$x\!\in\! X_{l+1}$, we denote by $x^-$ the unique element $X_l$ for which there is 
$\varepsilon\!\in\! 2^{\leq 1}$ with ${(x^-)\varepsilon\! =\! x}$. We also set $q_y\! :=\! p^{l+1}_{t_r1,x_{C(y)}}$ if $y\!\in\! C(t_r1)$. If $x\!\notin\! C(t_r1)$, or if $x\!\notin\! q_x$ and ${x^-\!\not=\! t_r}$, then we set 
$V_x\! :=\! U^l_{x^-}$. We also set $V_{t_r1}\! :=\! O_{t_r1}$ and $V_{t_r0}\! :=\! O_{t_r0}$. If now 
$x\!\in\! C(t_r1)$, $x\!\in\! q_x\!\setminus\!\{ t_r1\}$ and $x\! =\! q_x(i)$ with $i\!\geq\! 1$, then we define 
$V_x$ by induction on $i$. We set $V_x\! :=\! f_{u_l((q_x(i-1))^-)}[V_{q_x(i-1)}]$.

\vfill\eject

 This is possible since, inductively on $k\! <\!\vert q_x\vert$, $V_{q_x(k)}$ is defined, nonempty and contained in $W_{(q_x(k))^-}$, since $\Big(\big( q_x(k)\big)^-,\big( q_x(k\! +\! 1)\big)^-\Big)\!\in\! A_l$, so that 
$f_{u_l((q_x(k))^-)}[W_{(q_x(k))^-}]\! =\! W_{(q_x(k+1))^-}$.\bigskip

 Note that $\big( u,(V_x)_{x\in X_{l+1}}\big)\!\in\! U_{\cal T}$. Indeed, let $(y,x)\!\in\! A_{l+1}$. If 
$y\!\notin\! C(t_r1)$, then $y^-\!\not=\! t_r$, so that $(y^-,x^-)\!\in\! A_l$. Moreover, $x\!\notin\! C(t_r1)$ and $V_y\! =\! U^l_{y^-}$, $V_x\! =\! U^l_{x^-}$. Thus $V_y\!\subseteq\! D_{\phi (\varphi_{l-1}(y^-,x^-))}$ and 
$V_x\!\subseteq\! f_{\phi (\varphi_{l-1}(y^-,x^-))}[V_y]$, by the induction assumption. It remains to note that 
$$u(y)\! =\!\phi\big(\varphi_l(y,x)\big)\! =\!\phi\big(\varphi_{l-1}(y^-,x^-)\big)$$ 
to see that $V_y\!\subseteq\! D_{u(y)}$ and $V_x\!\subseteq\! f_{u(y)}[V_y]$. If $y\!\in\! C(t_r1)$ and 
$y\!\in\! q_y$, then $x\!\in\! q_y\!\setminus\!\{ t_r1\}$. By definition, $V_x\! =\! f_{u_l(y^-)}[V_y]$. Note that 
$u(y)\! =\!\phi\big(\varphi_l(y,x)\big)\! =\!\phi\big(\varphi_{l-1}(y^-,x^-)\big)\! =\! u_l(y^-)$. If $y\! =\! t_r0$, then 
$x\! =\! t_r1$, $V_{t_r0}\! =\! O_{t_r0}\!\subseteq\! D_{\phi (r)}\! =\! D_{u(t_r0)}$ and 
$V_{t_r1}\! =\! O_{t_r1}\!\subseteq\! f_{u(t_r0)}[V_{t_r0}]$. Otherwise, $V_y\! =\! U^l_{y^-}$ and 
$x\!\in\! C(t_r1)$. If $x\!\notin\! q_y$ and $x^-\!\not=\! t_r$, then we argue as in the case 
$y\!\notin\! C(t_r1)$. If $x\!\in\! q_y$ or $x^-\! =\! t_r$, then 
$V_x\!\subseteq\! W_{x^-}\! =\! f_{u_l(y^-)}[W_{y^-}]\!\subseteq\! f_{u_l(y^-)}[U^l_{y^-}]\! =\! f_{u(y)}[V_y]$.\bigskip 
 
 Lemma \ref{disj} provides a sequence $(U^{l+1}_x)_{x\in X_{l+1}}$ of pairwise disjoint nonempty clopen subsets of $Z$ with diameter at most $2^{-l-1}$ such that $U^{l+1}_y\!\subseteq\! D_{u(y)}$ if 
$y\!\notin\!\mbox{max}_{X_{l+1}}$ and 
$$U^{l+1}_x\!\subseteq\! V_x\cap\bigcap_{y\in\mbox{Pred}(x)}~f_{u(y)}[U^{l+1}_y]$$ 
if $x\!\in\! X_{l+1}$. Note that $U^{l+1}_x\!\subseteq\! V_x\!\subseteq\! U^l_{x^-}$, by definition of the $W_x$'s, which shows that (1) is satisfied. The construction of the $U^{l+1}_x$'s shows that (2)-(4) are satisfied. For (5), let $(y,x)\!\in\! B_{l+1}$, and $m\! <\! j$ with 
$\varphi_l\big( p(j\! -\! 1),p(j)\big)\! =\!\varphi_l(y,x)\! =\!\mbox{min}_{i<j}~\varphi_l\big( p(i),p(i\! +\! 1)\big)$. If $p(q)\!\not=\! t_r1$ for each $0\! <\! q\!\leq\! j$, then we are done, by induction assumption. If $p(q)\! =\! t_r1$ for some $0\! <\! q\!\leq\! j$, then $q\! <\! j$ or $j\!=\! 1$, and we may assume that $q\! <\! j$. It remains to note that $U^{l+1}_{p(q-\varepsilon )}\!\subseteq\! U^l_{(p(q))^-}\!\subseteq\! D_{\phi (\varphi_l(y,x))}$ for each 
$\varepsilon\!\in\! 2$, by the induction assumption.\hfill{$\square$}\bigskip

\noindent\bf Proof of Theorem \ref{cormain}.\rm\ We apply Lemma \ref{Bas0+} to $\big( X,(f_n)\big)$, 
$P\! :=\! X$ and $S_n\! :=\! D_n$. This is possible by Corollary \ref{D2+}. Lemma \ref{Bas0+} provides a Borel subset $S$ of $X$, a finer topology $\tau$ on $S$, and a sequence $(C_n)$ of clopen subsets of 
$Y\! :=\! (S,\tau )$ such that $\big( Y,({f_n}_{\vert C_n})\big)$ is a strongly complex situation and 
$C_n\!\subseteq\! S\cap g_n^{-1}(S)$. By Theorem \ref{newmain}, 
$(2^\omega,\mathbb{G}_1)\preceq^{\mbox{inj}}_c\big( Y,\bigcup_{n\in\omega}~\mbox{Graph}({f_n}_{\vert C_n})\big)$. As $\big( Y,\bigcup_{n\in\omega}~\mbox{Graph}({f_n}_{\vert C_n})\big)\preceq^{\mbox{inj}}_c(X,A^f)$, we are done.\hfill{$\square$}\bigskip

\noindent\bf Proof of Theorem \ref{G1}.\rm\ By Corollary \ref{D2+}, $\mathbb{G}_1$ is a $\boratwo$ digraph of uncountable Borel chromatic number, and is in particular analytic. Theorem \ref{cormain} shows that if $S$ is a Borel subset of $2^\omega$, $\tau$ is a finer topology on $S$, and $(C_n)$ is a sequence of clopen subsets of $Y\! :=\! (S,\tau )$ such that $\big( Y,({g_n}_{\vert C_n})\big)$ is a strongly complex situation, then $(2^\omega ,\mathbb{G}_1)\preceq^{\mbox{inj}}_c
\big( Y,\bigcup_{n\in\omega}~\mbox{Graph}({g_n}_{\vert C_n})\big)$. We apply Corollary \ref{charmin+} and the remark after it to get the minimality of $\mathbb{G}_1$. This implies the minimality of $\mathbb{G}_1^{-1}$. We saw in the introduction that $\mathbb{G}_0$ is also minimal. It remains to apply Lemma \ref{propf++} and Corollary \ref{th12}.\hfill{$\square$}

\section{$\!\!\!\!\!\!$ The structure of $\preceq^{\mbox{inj}}_{\mathfrak{C}}$}

\begin{thm} \label{incomp} $(2^\omega ,\mathbb{G}_L)_{L\in\omega}$ is a sequence made of $\boratwo$ digraphs pairwise 
$\preceq^{\mbox{inj}}_B$-incompatible among analytic digraphs of uncountable Borel chromatic number on Polish spaces. In particular, the sequence $(2^\omega ,\mathbb{G}_L)_{L\in\omega}$ is a $\preceq^{\mbox{inj}}_B$-antichain. Moreover, if $L\!\not=\! M$ are natural numbers, $s\!\in\! 2^{<\omega}$ and $G$ is a dense $G_\delta$ subset of $2^\omega$, then 
$\big(N_s\cap G,\mathbb{G}_L\cap (N_s\cap G)^2\big)\not\preceq^{\mbox{inj}}_c(2^\omega ,\mathbb{G}_M)$.\end{thm}

\noindent {\bf Proof.} By Lemma \ref{propf++}.(a), $\big( 2^\omega ,({g^L_n}_{\vert\mathbb{D}^L_n})\big)$ is a complex situation if $L\!\geq\! 1$. By Theorem \ref{G0} and Corollary \ref{D2+}, $\mathbb{G}_L$ is a $\boratwo$ digraph on $2^\omega$ of uncountable Borel chromatic number.\bigskip

\noindent {\bf Claim.} \it\ Let $L\! <\! M$ be natural numbers, and $\big( X,(f_n)\big)$ be a complex situation satisfying the Conditions (a), (b) of a strongly complex situation. Then $(X,A^f)\preceq_c^{\mbox{inj}}(2^\omega ,\mathbb{G}_L),(2^\omega ,\mathbb{G}_M)$ cannot hold simultaneously.\rm\bigskip

 Indeed, we argue by contradiction, which gives witnesses $u,u'$. We set\bigskip

\leftline{$R\! :=\!\{ (n,r,r',x)\!\in\!\omega^3\!\times\! X\mid x\!\in\! D_n\wedge u(x)\!\in\!\mathbb{D}_r^L\wedge 
u'(x)\!\in\!\mathbb{D}_{r'}^M\wedge u\big( f_n(x)\big)\! =\! g^L_r\big( u(x)\big)~\wedge$}\smallskip

\rightline{$u'\big( f_n(x)\big)\! =\! g^M_{r'}\big( u'(x)\big)\} .$}\bigskip 

\noindent Note that $R$ is closed, by continuity. If $n\!\in\!\omega$, $U$ is a nonempty clopen subset of $D_n$ and $x\!\in\! U$, then there are 
$r,r'\!\in\!\omega$ such that $(n,r,r',x)\!\in\! R$. By Baire's theorem, we can find a nonempty clopen subset $C$ of $U$ and $r,r'\!\in\!\omega$ such that 
$(n,r,r',x)\!\in\! R$ if $x\!\in\! C$.\bigskip

 We inductively construct a strictly increasing sequence $(n_k)_{k\in\omega}\!\in\!\omega^\omega$ as follows. We first apply the previous point to $n_0\! :=\! 0$ and $U\! :=\! D_0$, which provides $C_0\!\subseteq\! D_{n_0}$ and $r_0,r'_0\!\in\!\omega$ such that $(n_0,r_0,r'_0,x)\!\in\! R$ if $x\!\in\! C_0$. As 
$\big( X,(f_n)\big)$ is a strongly complex situation, we can find $n_{k+1}\! >\! n_k$ such that $\mbox{Graph}(f_{n_{k+1}})\cap C_k^2\!\not=\!\emptyset$. We apply the previous point to $n_{k+1}$ and $U\!\subseteq\! C_k\cap f_{n_{k+1}}^{-1}(C_k)$ with diameter at most $2^{-k}$, which provides 
$C_{k+1}\!\subseteq\! U$ and $r_{k+1},r'_{k+1}\!\in\!\omega$ such that $(n_{k+1},r_{k+1},r'_{k+1},x)$ is in $R$ if $x\!\in\! C_{k+1}$.\bigskip 

 As $(C_k)$ is a decreasing sequence of nonempty closed subsets of $X$ whose diameters tend to zero, there is $x\!\in\!\bigcap_{k\in\omega}~C_k$. Note that 
$(n_k,r_k,r'_k,x)\!\in\! R$ for each $k$. Let us prove that $(r_k)$ is unbounded. We argue by contradiction, which gives $r\!\in\!\omega$ and 
$I\!\subseteq\!\omega$ infinite such that $r_k\! =\! r$ if $k\!\in\! I$. By continuity, we get $u(x)\! =\! g^L_r\big( u(x)\big)$ since 
$\mbox{lim}_{k\rightarrow\infty ,k\in I}~f_{n_k}(x)\! =\! x$, which contradicts the fact that $g^L_r$ is fixed point free. So, extracting a subsequence if necessary, we may assume that $(r_k)$ and $(r'_k)$ are strictly increasing.\bigskip

 As $x\!\in\! C_1$ and $f_{n_1}(x)\!\in\! C_0$, $u\Big( f_{n_0}\big( f_{n_1}(x)\big)\Big)\! =\! 
g^L_{r_0}\Big( u\big( f_{n_1}(x)\big)\Big)\! =\! g^L_{r_0}\Big( g^L_{r_1}\big( u(x)\big)\Big)$. More generally, 
$u\big( f_{n_s}(x)\big)\! =\! g^L_{r_s}\big( u(x)\big)$ and $u'\big( f_{n_s}(x)\big)\! =\! g^M_{r'_s}\big( u'(x)\big)$ if 
$s\!\in\!\omega^{<\omega}$ is strictly increasing.\bigskip

 Assume first that $L\!\geq\! 1$. We choose $s\!\in\!\omega^{L+1}$ strictly increasing. By Lemma \ref{propf++}.(c), we get 
$g^L_{r_s}\big( u(x)\big)\! =\! g^L_{r_{s^-}}\big( u(x)\big)$, so that $g^M_{r'_s}\big( u'(x)\big)\! =\! g^M_{r'_{s^-}}\big( u'(x)\big)$ by injectivity of $u$, which contradicts Lemma \ref{propf++}.(b) since $u'(x)\!\in\!\mathbb{D}_{r_{s(\vert s\vert -1)}}^M$. If 
$L\! =\! 0$, we argue simililarly, using the fact that the $g^0_r$'s are injective and fixed point free.\hfill{$\diamond$}

\vfill\eject

Let $L\! <\! M$. We argue by contradiction, which gives a Polish space $Y$ and an analytic digraph $B$ on $Y$ of uncountable Borel chromatic number such that $(Y,B)\preceq^{\mbox{inj}}_B(2^\omega ,\mathbb{G}_L),(2^\omega ,\mathbb{G}_M)$. Theorem \ref{fav} gives a complex situation $\big( X,(f_n)\big)$ such that $(X,A^f)\preceq_c^{\mbox{inj}}(Y,B)$ or $(X,A^f)\preceq_c^{\mbox{inj}}(Y,B^{-1})$.\bigskip

 Let us prove that $(X,A^f)\preceq_c^{\mbox{inj}}(Y,B)$. We argue by contradiction, so that 
$$(X,A^f)\preceq_c^{\mbox{inj}}(Y,B^{-1}).$$ 
As $\mathbb{G}_L$ has countable vertical sections, $B$ too, so that $B^{-1}$ and $A^f$ have countable horizontal sections. Thus $A^f$ is locally countable. By Corollary \ref{D2+}, $A^f$ is a $\boratwo$ digraph of uncountable Borel chromatic number. The discussion after Theorem \ref{G0} shows that 
$(2^\omega ,\mathbb{G}_0)\preceq_c^{\mbox{inj}}(X,A^f)$, and also 
$$(2^\omega ,\mathbb{G}_0)\preceq_c^{\mbox{inj}}(2^\omega ,\mathbb{G}_0^{-1})\preceq_c^{\mbox{inj}}(Y,B)
\preceq^{\mbox{inj}}_B(2^\omega ,\mathbb{G}_L),(2^\omega ,\mathbb{G}_M).$$ 
Let $H$ be a dense $G_\delta$ subset of $2^\omega$ such that 
$(H,\mathbb{G}_0\cap H^2)\preceq_c^{\mbox{inj}}(2^\omega ,\mathbb{G}_L),(2^\omega ,\mathbb{G}_M)$. The proof of Lemma \ref{suffnpc}.(b) shows that 
$\mathbb{G}_0\cap H^2$ has uncountable Borel chromatic number. The discussion after Theorem \ref{G0} shows that $(2^\omega ,\mathbb{G}_0)\preceq_c^{\mbox{inj}}(H,\mathbb{G}_0\cap H^2)$, and thus 
$(2^\omega ,\mathbb{G}_0)\preceq_c^{\mbox{inj}}(2^\omega ,\mathbb{G}_L),(2^\omega ,\mathbb{G}_M)$, which contradicts the claim.\bigskip

 This shows that $(X,A^f)\preceq_B^{\mbox{inj}}(2^\omega ,\mathbb{G}_L),(2^\omega ,\mathbb{G}_M)$. So we may assume that $Y\! =\! X$ and $B\! =\! A^f$. Let $P$ be a dense $G_\delta$ subset of $X$ such that 
$(P,A^f\cap P^2)\preceq_c^{\mbox{inj}}(2^\omega ,\mathbb{G}_L),(2^\omega ,\mathbb{G}_M)$. The proof of Lemma \ref{suffnpc}.(b) shows that $A^f\cap P^2$ has uncountable Borel chromatic number. The previous point shows that $(2^\omega ,\mathbb{G}_0)\not\preceq_c^{\mbox{inj}}(X,A^f)$. By Lemma \ref{Bas0+}, we may assume that $\big( X,(f_n)\big)$ is a strongly complex situation and 
$(X,A^f)\preceq_c^{\mbox{inj}}(2^\omega ,\mathbb{G}_L),(2^\omega ,\mathbb{G}_M)$, which contradicts the claim.\bigskip

 The last assertion comes from the fact that $\mathbb{G}_L\cap (N_s\cap G)^2$ has uncountable Borel chromatic number, by the proof of Lemma 
\ref{suffnpc}.(b).\hfill{$\square$}\bigskip

\noindent {\bf Proof of Theorem \ref{infi}.} We argue by contradiction, which gives a natural number $N$ and a basis 
$(X_i,A_i)_{i<N}$. By Theorem \ref{incomp}, we can find $i\! <\! N$ and $L\!\not=\! M$ with 
$(X_i,A_i)\preceq^{\mbox{inj}}_B(2^\omega ,\mathbb{G}_L),(2^\omega ,\mathbb{G}_M)$. This contradicts Theorem \ref{incomp}.\hfill{$\square$}\bigskip

 We now prove Theorem \ref{biganti}.\bigskip

\noindent {\bf Notation.} We set, for each $n\!\in\!\omega$, $\mathbb{Q}_n\! :=\!\oplus_{L\leq n}~2^\omega$ and  
$$\mathbb{H}_n\! :=\!\big\{\big( (L,\gamma ),(M,\delta )\big)\!\in\!\mathbb{Q}_n^2\mid L\! =\! M\wedge 
(\gamma ,\delta )\!\in\!\mathbb{G}_L\big\} .$$ 

\begin{thm} Let ${\mathfrak C}\!\in\!\{ c,B\}$. Then $(\mathbb{Q}_n,\mathbb{H}_n)_{n\in\omega}$ is a 
$\preceq^{\mbox{inj}}_{\mathfrak{C}}$-strictly increasing chain made of $\boratwo$ digraphs of uncountable Borel chromatic number.\end{thm}

\noindent {\bf Proof.} By Theorem \ref{incomp}, $\mathbb{H}_n$ is a $\boratwo$ digraph of uncountable Borel chromatic number on $\mathbb{Q}_n$. The identity map shows that our sequence is increasing. We argue by contradiction, which gives $n\!\in\!\omega$ such that 
$(\mathbb{Q}_{n+1},\mathbb{H}_{n+1})\preceq^{\mbox{inj}}_B(\mathbb{Q}_n,\mathbb{G}_n)$ with witness $u$. The map associating $u_0(n,\gamma )$ to $\gamma\!\in\! 2^\omega$ is Borel, which gives $L\! <\! n$, 
$s\!\in\! 2^{<\omega}$ and a dense $G_\delta$ subset $G$ of $2^\omega$ such that $u_0(L,\gamma )\! =\! L$ if 
$\gamma$ is in $N_s\cap G$. 

\vfill\eject

 The map $v$ associating $u_1(n,\gamma )$ to $\gamma\!\in\! N_s\cap G$ is injective and Borel. We can restrict $G$, so that we may assume that $v$ is continuous. Thus 
$\big( N_s\cap G,\mathbb{G}_n\cap (N_s\cap G)^2\big)\preceq^{\mbox{inj}}_c(2^\omega ,\mathbb{G}_L)$ with witness 
$v$, which contradicts Theorem \ref{incomp} since $n\!\not=\! L$.\hfill{$\square$}\bigskip

\noindent {\bf Notation.} Let $(p_n)_{n\in\omega}$ be the sequence of prime numbers. We define, for each $\alpha\!\in\! 2^\omega$, 
$E_\alpha\!\subseteq\!\omega$ by $E_\alpha\! :=\!\{ p_0^{\alpha (0)+1}\cdots p_n^{\alpha (n)+1}\mid n\!\in\!\omega\}$. Then we set $\mathbb{P}_\alpha\! :=\!\oplus_{L\in E_\alpha}~2^\omega$ and  
$$\mathbb{G}_\alpha\! :=\!\big\{\big( (L,\gamma ),(M,\delta )\big)\!\in\!\mathbb{P}_\alpha^2\mid L\! =\! M\wedge (\gamma ,\delta )\!\in\!\mathbb{G}_L\big\} .$$  

\begin{thm} \label{antichain} $(\mathbb{P}_\alpha ,\mathbb{G}_\alpha )_{\alpha\in 2^\omega}$ is a 
$\preceq^{\mbox{inj}}_B$-antichain made of $\boratwo$ digraphs of uncountable Borel chromatic number.\end{thm}

\noindent {\bf Proof.} By Theorem \ref{incomp}, $\mathbb{G}_\alpha$ is a $\boratwo$ digraph of uncountable Borel chromatic number on $\mathbb{P}_\alpha$. We argue by contradiction, which gives $\alpha\!\not=\!\beta$ such that 
$(\mathbb{P}_\alpha ,\mathbb{G}_\alpha )\preceq^{\mbox{inj}}_B(\mathbb{P}_\beta ,\mathbb{G}_\beta )$ with witness 
$u$. Note that $E_\alpha\cap E_\beta$ is finite, which gives $L\!\in\! E_\alpha\!\setminus\! E_\beta$. The map associating $u_0(L,\gamma )$ to $\gamma\!\in\! 2^\omega$ is Borel, which gives $M\!\in\! E_\beta$, 
$s\!\in\! 2^{<\omega}$ and a dense $G_\delta$ subset $G$ of $2^\omega$ such that 
$u_0(L,\gamma )\! =\! M$ if $\gamma$ is in $N_s\cap G$. The map $v$ associating $u_1(L,\gamma )$ to 
$\gamma\!\in\! N_s\cap G$ is injective and Borel. We can moreover restrict $G$, so that we may and will assume that $v$ is continuous. Thus $v$ is a witness for the fact that 
$\big( N_s\cap G,\mathbb{G}_L\cap (N_s\cap G)^2\big)\preceq^{\mbox{inj}}_c
(2^\omega ,\mathbb{G}_M)$, which contradicts Theorem \ref{incomp} since $L\!\not=\! M$.\hfill{$\square$}\bigskip

 Our main results also hold for graphs.

\begin{thm} \label{graphs} Let ${\mathfrak C}\!\in\!\{ c,B\}$.\smallskip

(a) There is a $\preceq^{\mbox{inj}}_{\mathfrak C}$-antichain $\{ s(\mathbb{G}_0),s(\mathbb{G}_1)\}$ made of graphs $\preceq^{\mbox{inj}}_{\mathfrak C}$-minimal among analytic graphs of uncountable Borel chromatic number.\smallskip

(b) There is a $\preceq^{\mbox{inj}}_{\mathfrak{C}}$-antichain of size $2^{\aleph_0}$ made of 
$\boratwo$ graphs of uncountable Borel chromatic number.\smallskip

(c) Any $\preceq^{\mbox{inj}}_{\mathfrak{C}}$-basis for the class of analytic graphs of uncountable Borel chromatic number on Polish spaces is infinite.\end{thm}

\noindent {\bf Proof.} We first prove the following.\bigskip

\noindent {\bf Claim.}\it\ $\big(2^\omega ,s(\mathbb{G}_L)\big)_{L\in\omega}$ is a sequence made of $\boratwo$ graphs pairwise $\preceq^{\mbox{inj}}_B$-incompatible among analytic graphs of uncountable Borel chromatic number on Polish spaces.\rm\bigskip

 Indeed, we argue by contradiction, which gives $L\!\not=\! M$ and an analytic graph $G$ of uncountable Borel chromatic number on a Polish space $X$ such that 
$(X,G)\preceq^{\mbox{inj}}_B\big(2^\omega ,s(\mathbb{G}_L)\big),\big(2^\omega ,s(\mathbb{G}_M)\big)$, with witnesses $u$, $u'$ respectively. Note that 
$$\begin{array}{ll}
G\! =\!\!\!\!\!
& \big( G\cap (u\!\times\! u)^{-1}(\mathbb{G}_L)\cap (u'\!\times\! u')^{-1}(\mathbb{G}_M)\big)\cup
\big( G\cap (u\!\times\! u)^{-1}(\mathbb{G}_L)\cap (u'\!\times\! u')^{-1}(\mathbb{G}^{-1}_M)\big)\ \cup\cr
& \big( G\cap (u\!\times\! u)^{-1}(\mathbb{G}^{-1}_L)\cap (u'\!\times\! u')^{-1}(\mathbb{G}_M)\big)\cup
\big( G\cap (u\!\times\! u)^{-1}(\mathbb{G}^{-1}_L)\cap (u'\!\times\! u')^{-1}(\mathbb{G}^{-1}_M)\big) .
\end{array}$$
The second and the third of these subgraphs are locally countable. If one of them has uncountable Borel chromatic number, then it is above $\mathbb{G}_0$ by the discussion after Theorem \ref{G0}. By Theorem \ref{incomp}, we must have $L\! =\! 0\! =\! M$, which is absurd. Thus the first or the fourth of these subgraphs has uncountable Borel chromatic number, which contradicts the incompatibility of $\mathbb{G}_L$ and $\mathbb{G}_M$.\hfill{$\diamond$}

\vfill\eject

\noindent (a) By the claim, $s(\mathbb{G}_0)$ and $s(\mathbb{G}_1)$ are incompatible, and thus incomparable. Assume that we can find $L\!\in\! 2$ and an analytic graph $G$ of uncountable Borel chromatic number on a Polish space $X$ such that $(X,G)\preceq^{\mbox{inj}}_{\mathfrak{C}}\big(2^\omega ,s(\mathbb{G}_L)\big)$, with witness $u$. Note that $$G\! =\!\big( G\cap (u\!\times\! u)^{-1}(\mathbb{G}_L)\big)\cup
\big( G\cap (u\!\times\! u)^{-1}(\mathbb{G}^{-1}_L)\big) .$$ 
One of these subgraphs has uncountable Borel chromatic number. Assume for example that it is the first one. Then the minimality of $\mathbb{G}_L$ shows that $(2^\omega ,\mathbb{G}_L)\preceq^{\mbox{inj}}_{\mathfrak{C}}
\Big( 2^\omega ,\big( G\cap (u\!\times\! u)^{-1}(\mathbb{G}_L)\big)\Big)$, and thus 
$\big(2^\omega ,s(\mathbb{G}_L)\big)\preceq^{\mbox{inj}}_{\mathfrak{C}}(X,G)$.\bigskip

\noindent (b) We essentially argue as in the proof of Theorem \ref{antichain}. By the claim, $s(\mathbb{G}_\alpha )$ is a $\boratwo$ graph of uncountable Borel chromatic number on $\mathbb{P}_\alpha$. We argue by contradiction, which gives $\alpha\!\not=\!\beta$ such that $\big(\mathbb{P}_\alpha ,s(\mathbb{G}_\alpha )\big)
\preceq^{\mbox{inj}}_B\big(\mathbb{P}_\beta ,s(\mathbb{G}_\beta )\Big)$ with witness $u$. Note that $E_\alpha\cap E_\beta$ is finite, which gives $L\!\in\! E_\alpha\!\setminus\! E_\beta$. The map associating $u_0(L,\gamma )$ to $\gamma\!\in\! 2^\omega$ is Borel, which gives $M\!\in\! E_\beta$, 
$s\!\in\! 2^{<\omega}$ and a dense $G_\delta$ subset $G$ of $2^\omega$ such that 
$u_0(L,\gamma )\! =\! M$ if $\gamma$ is in $N_s\cap G$. The map $v$ associating $u_1(L,\gamma )$ to 
$\gamma\!\in\! N_s\cap G$ is injective and Borel. We can moreover restrict $G$, so that we may and will assume that $v$ is continuous. Thus $v$ is a witness for the fact that 
$\big( N_s\cap G,s(\mathbb{G}_L)\cap (N_s\cap G)^2\big)\preceq^{\mbox{inj}}_c
\big(2^\omega ,s(\mathbb{G}_M)\big)$, which contradicts the claim \ref{incomp} since $L\!\not=\! M$.\bigskip

\noindent (c) We argue as in the proof of Theorem \ref{infi}, using the claim.\hfill{$\square$}

\section{$\!\!\!\!\!\!$ References}\indent
 
\noindent [K]\ \ A. S. Kechris,~\it Classical Descriptive Set Theory,~\rm 
Springer-Verlag, 1995

\noindent [K-S-T]\ \ A. S. Kechris, S. Solecki and S. Todor\v cevi\'c,~\it Borel chromatic numbers,\ \rm 
Adv. Math.~141 (1999), 1-44

\noindent [K-Ma]\ \ A. S. Kechris, and A. S. Marks, Descriptive graph combinatorics,\ \it preprint (see the url http://www.math.caltech.edu/$\sim$kechris/papers/combinatorics18.pdf)\ \rm

\noindent [L1]\ \ D. Lecomte,~\it Classes de Wadge potentielles et th\'eor\`emes d'uniformisation 
partielle,~\rm Fund. Math.~143 (1993), 231-258

\noindent [L2]\ \ D. Lecomte, Uniformisations partielles et crit\`eres \`a la Hurewicz dans le plan,\ \it Trans. Amer. Math. Soc.\rm\ 347, 11 (1995), 4433-4460

\noindent [L3]\ \ D. Lecomte,~\it Tests \`a la Hurewicz dans le plan,~\rm Fund. Math.~156 (1998), 
131-165

\noindent [L4]\ \ D. Lecomte, On minimal non potentially closed subsets of the plane,\ \it Topology Appl.\rm\ 154, 1 (2007), 241-262

\noindent [L-Mi]\ \ D. Lecomte and B. D. Miller, Basis theorems for non-potentially closed sets and graphs of uncountable Borel chromatic number,\ \it J. Math. Log.\ \rm 8 (2008), 1-42

\noindent [Lo]\ \ A. Louveau,~\it Ensembles analytiques et bor\'eliens dans les espaces produit,~\rm Ast\'erisque (S. M. F.) 78 (1980)

\noindent [Mi]\ \ B. D. Miller, The graph-theoretic approach to descriptive set theory,\ \it Bull. Symbolic Logic\ \rm 18, 4 (2012), 554-575

\noindent [M]\ \ Y. N. Moschovakis,~\it Descriptive set theory,~\rm North-Holland, 1980

\end{document}